# NASH EQUILIBRIUM AND ROBUST STABILITY IN DYNAMIC GAMES: A SMALL-GAIN PERSPECTIVE

**Iasson Karafyllis[*], Zhong-Ping Jiang[**] and George Athanasiou[***]**

[*]Department of Environmental Engineering,
Technical University of Crete, 73100, Chania, Greece
email: ikarafyl@enveng.tuc.gr

[**]Department of Electrical and Computer Engineering,
Polytechnic Institute of New York University,
Six Metrotech Center, Brooklyn, NY 11201, U.S.A.
email: zjiang@control.poly.edu

[***]TT Hellenic Postbank, Financial Services,
2-6 Pesmazoglou str., 101 75, Athens, Greece
email: g.athanasiou@ttbank.gr

**Abstract**

This paper develops a novel methodology to study robust stability properties of Nash equilibrium points in dynamic games. Small-gain techniques in modern mathematical control theory are used for the first time to derive conditions guaranteeing uniqueness and global asymptotic stability of Nash equilibrium point for economic models described by functional difference equations. Specification to a Cournot oligopoly game is studied in detail to demonstrate the power of the proposed methodology.

## 1. Introduction

Dynamical game-theoretical models have inherent uncertainty in many aspects. The uncertainty is related strongly to a number of open questions which cannot be answered a priori:

1. Should the models be formulated in continuous time or discrete time?
   The answer to the above question is crucial: models in discrete time will be described by difference equations (see [1,2,3,5,6,11,21,33]) while models in continuous time are generally described by differential equations (with or without delays; see [4,31]). The answer to the above question has significant consequences: the perception of time for each player in a dynamic game-theoretical model affects her behavior.

2. What are the expectation rules that a player has for the other players?
   Again the answer to the above question is crucial: the behavior of a player will heavily rely on the expectations for the actions of the rest players. There is a large economic literature on the effect of expectation rules (e.g., naïve, backward-looking, rational expectations, see [3,5,10,11] and references therein). Moreover, if expectation rules are using delay terms then the consequences on stability can be important (see [4,20]).

3. What are the values of the various constants involved in a dynamic game-theoretical model?
   In many dynamic games, the rate of change of the action of one player is assumed to be proportional to either the deviation of the action from the best reply (see for example [6,11]) or the gradient of the payoff function (see for example [6,31]). The value of the proportionality constant cannot be known a priori.

Therefore, the answers to important questions such as the existence of a Nash equilibrium point, the uniqueness of a Nash equilibrium point and its stability properties are usually related to the specific assumptions made in order to cope with the uncertainty. Consequently, the following question arises:

"Can we extract robust information from an uncertain nonlinear economic model,
which will hold no matter what the uncertainty is?"

The present work answers it affirmatively. In some cases, we can even show the existence of a Nash equilibrium point, its uniqueness and its global asymptotic stability properties for all possible uncertainties. In order to be able to do this we propose the following methodology:

- First, we formulate our models in continuous time by means of Functional Difference Equations (see [9,15-19,32]). By doing so we convert a finite-dimensional problem to an infinite-dimensional problem, which seems to be a clear disadvantage at first sight. However, in this way we can obtain all features of continuous time and discrete time models. Indeed, we will show that many models appeared in the literature can be considered as special cases of our proposed model.

- Second, we do not assume a specific expectation rule: instead, we will only assume that the expectation is consistent with the history of the game (consistent backward looking expectation; see Definition 2.1 below).

- In order to be able to extract important information from the uncertain model we use advanced stability methods. Indeed, by applying small-gain analysis (see [17,18,19]), we can guarantee that the Nash equilibrium point is unique and globally asymptotically stable (see Theorem 3.1 and Theorem 4.2 below).

To our knowledge, this is the first time that such results are presented for dynamical game-theoretical nonlinear models. The only other work, which we have found and can address such questions, is [21]: our results generalize the results in [21]. Moreover, the results of [21] are applied in a discrete-time framework and cannot be used for the analysis of models in continuous time. As a byproduct of our work, we will also give conditions for uniqueness of a fixed point (see Corollary 3.2 and Corollary 3.3 below), which can be used in conjunction with classical fixed-point theorems and are different from other uniqueness conditions in the literature (see [8]).

It should be noticed that the stability/uniqueness conditions obtained by the proposed methodology will be more demanding conditions than the ones which can be obtained from the study of a specific model (with specific expectation rules, specific values for the constants involved in the model and with a specific perception of time). However, this is expected since the stability/uniqueness conditions obtained by the proposed methodology are sufficient conditions for global asymptotic stability for an uncertain model, which contains many other models as special cases. To this end stability analysis by means of nonlinear small-gain theorems is utilized. Small-Gain results have been used frequently in stability studies (see [7,12,13,14,17,18,19]) and are based on variations of the Input-to-State Stability property introduced by E. D. Sontag in [26] and the Input-to-Output Stability property (see [14,16,27,28]).

The structure of the paper is as follows: in Section 2, we apply the above described methodology to the Cournot dynamic oligopoly problem. There is a vast literature on this well-studied problem (see for instance [1,3,4,5,6,11,29,31,33]). For this specific problem, we describe in detail our proposed methodology and we show how we can obtain results on the stability properties of the Cournot equilibrium, which do not depend on the form of the uncertainty. The presentation of the special case of the Cournot game before the general case was preferred for tutorial purposes: all issues arising in the general case are present in the Cournot game. In Section 3, we proceed to the more general case of dynamic strategic games and in Section 4 we discuss the problem of accommodating the rational expectations. Our concluding remarks are given in Section 5. Finally, in the Appendix, we give the proofs of certain results of this work.

**Notations** Throughout this paper we adopt the following notations:

∗ For a vector $x \in \Re^n$ we denote by $|x|$ its usual Euclidean norm.

∗ $\Re^+$ denotes the set of non-negative real numbers. For every $t \in \Re^+$, $[t]$ denotes the integer part of $t$, i.e., the largest integer being less than or equal to $t$.

∗ We say that a non-decreasing continuous function $\gamma : \Re^+ \to \Re^+$ is of class $N$ if $\gamma(0) = 0$.

∗ Let $x : I \to \Re^n$ with $[a,b] \subseteq I$ and $\sup_{\tau \in I} |x(\tau)| < +\infty$. We denote by $\|x\|_{[a,b]} = \sup_{a \leq \tau \leq b} |x(\tau)|$.

∗ Let $U \subseteq \Re^n$ be a closed convex set. By $\Pr_U(x)$ we denote the projection of $x \in \Re^n$ on $U \subseteq \Re^n$.

* The norm of a normed linear space $X$ will be denoted by $\| \ \|_X$. More specifically, in the present work $X$ will denote the normed linear space of bounded functions $x:[-T,0] \to \Re^n$ with norm $\|x\|_X = \sup_{-T \leq \tau \leq 0} |x(\tau)|$, for given $T \geq 0$. If $x:[-T,a] \to \Re^n$, where $a \geq 0$, is a bounded mapping then $x_t \in X$ with $t \in [0,a]$ is defined by $x_t = \{x(\tau): t-T \leq \tau \leq t\}$ as usually in systems with delays (see [9]).

* For a vector $q = (q_1,...,q_n) \in S_1 \times ... \times S_n$ we will use the notation (see [30])

$$q_{-i} = (q_1,...,q_{i-1},q_{i+1},...,q_n) \text{ for } 1 < i < n \text{ and } n \geq 3$$

$$q_{-1} = (q_2,...,q_n), \ q_{-n} = (q_1,...,q_{n-1}) \text{ for } n \geq 2$$

i.e., $q_{-i}$ is the vector of order $n-1$ after deleting the i-th component $q_i \in S_i$ of the vector $q = (q_1,...,q_n) \in S_1 \times ... \times S_n$.

## 2. Dynamic Cournot Oligopoly

We consider the case of Cournot oligopoly where $n$ players produce quantities of a single homogeneous product. The payoff function for each player is expressed by:

$$\pi_i = pq_i - c_i q_i - \frac{1}{2} K_i q_i^2, \ i = 1,...,n \tag{2.1}$$

where $K_i, c_i$, $i = 1,...,n$ are constants, $q_i \in [0, Q_i]$, $i = 1,...,n$ is the quantity of the commodity produced by the i-*th* player, $Q_i > 0$ is the maximum level of production of the product for the i-*th* player and $p \geq 0$ is the price of the commodity.

Assuming a linear demand function:

$$p = b\left(a - \sum_{i=1}^{n} q_i\right) \tag{2.2}$$

where $a, b > 0$ are constants satisfying $a \geq \sum_{i=1}^{n} Q_i$ and $b > -\frac{1}{2} \min_{i=1,...,n} K_i$, we obtain the best reply mapping for each one of the players:

$$q_i = f_i(q_{-i}) := \min\left\{Q_i, \max\left\{0, \frac{ab - c_i}{2b + K_i} - \frac{b}{2b + K_i} \sum_{j \neq i} q_j\right\}\right\}, \ i = 1,...,n \tag{2.3}$$

We define:

$$S := [0,Q_1] \times [0,Q_2] \times ... \times [0,Q_n] \subset \Re^n \tag{2.4}$$

$$q = (q_1,...,q_n) \in S \tag{2.5}$$

$$F(q) := \begin{bmatrix} f_1(q_{-1}) \\ \vdots \\ f_n(q_{-n}) \end{bmatrix} = \begin{bmatrix} \min\left\{Q_1, \max\left\{0, \frac{ab-c_1}{2b+K_1} - \frac{b}{2b+K_1} \sum_{j \neq 1} q_j\right\}\right\} \\ \vdots \\ \min\left\{Q_n, \max\left\{0, \frac{ab-c_n}{2b+K_n} - \frac{b}{2b+K_n} \sum_{j \neq n} q_j\right\}\right\} \end{bmatrix} \tag{2.6}$$

and we notice that the set $S \subset \Re^n$ as defined by (2.4) is compact and convex and that the map $F: S \to S$ as defined by (2.6) is continuous. Consequently, Brouwer's fixed point theorem guarantees the existence of at least one Nash equilibrium $q^* \in S$ with $q_i^* = f_i(q_{-i}^*) = \min\left\{Q_i, \max\left\{0, \frac{ab-c_i}{2b+K_i} - \frac{b}{2b+K_i}\sum_{j \neq i} q_j^*\right\}\right\}$ for $i = 1,...,n$.

Next we assume that the dynamics of the game are described in continuous time as follows:

- every player forms an expectation for the behavior of all other players at each time $t \geq 0$: the expectation of the i-*th* player for the production level of the j-*th* player at time $t \geq 0$ will be denoted by $q_{i,j}^{\exp}(t) \in [0, Q_j]$ ($j \neq i$, $i, j = 1,...,n$),
- every player determines her production level as a convex combination of a past production level and the best reply response based on the expectations for the behavior of all other players at each time $t \geq 0$, i.e.,

$$q_i(t) = \theta_i(t) \min\{Q_i, \max\{0, q_i(t - \tau_i(t))\}\} + (1 - \theta_i(t)) \min\left\{Q_i, \max\left\{0, \frac{ab-c_i}{2b+K_i} - \frac{b}{2b+K_i}\sum_{j \neq i} q_{i,j}^{\exp}(t)\right\}\right\}, \; i = 1,...,n$$

(2.7)

where $\theta_i : \Re^+ \to [0, \Theta]$, $\tau_i : \Re^+ \to [r, T]$, $i = 1,...,n$ are in general unknown functions, $0 \leq \Theta < 1$, $0 < r \leq T$ are constants (in general unknown).

The reader should notice that (2.7) is a model that evolves in continuous time, i.e., $t \in \Re^+$. If the expectation rules $q_{i,j}^{\exp}(t)$ ($j \neq i$, $i, j = 1,...,n$), and the functions $\theta_i : \Re^+ \to [0, \Theta]$, $\tau_i : \Re^+ \to [r, T]$ ($i = 1,...,n$) were known, we would have an accurate description of the dynamics of the Cournot oligopoly game. However, we will not assume exact knowledge of the expectation rules but a specific consistency condition. First we give the definition for a Consistent Backward-looking expectation with respect to the Nash equilibrium point $q^* \in S$.

**Definition 2.1:** *An expectation rule $q_{i,j}^{\exp}(t)$ (where $j \neq i$, $i, j = 1,...,n$) is called a Consistent Backward-looking expectation with respect to the Nash equilibrium point $q^* \in S$ if there exist constants $0 < r \leq T$ such that:*

$$\left|q_{i,j}^{\exp}(t) - q_j^*\right| \leq \sup_{t-T \leq \tau \leq t-r} \left|q_j(\tau) - q_j^*\right| = \left\|q_j - q_j^*\right\|_{[t-T, t-r]}, \text{ for all } t \geq 0 \quad (2.8)$$

In other words the consistency condition (2.8) recognizes that it is not logical for i-*th* player to expect that the production level of the j-*th* manufacturer will deviate from its equilibrium level more than the highest deviation she has experienced in the past. Next we present some examples of Consistent Backward-looking expectation rules:

1) $q_{i,j}^{\exp}(t) = a_{i,j}(t) \sum_{l=1}^{m} w_{i,j,l}(t) q_j(t - \tau_{i,j,l}(t)) + (1 - a_{i,j}(t)) q_j^*$, where $a_{i,j}(t) \in [0,1]$, $T \geq \tau_{i,j,l}(t) \geq r > 0$, $w_{i,j,l}(t) \geq 0$ with $1 = \sum_{l=1}^{m} w_{i,j,l}(t)$ for all $t \geq 0$ and $l = 1,...,m$. In discrete-time models the case $\tau_{i,j,l}(t) = t + l - [t]$, $a_{i,j}(t) \equiv 1$, $w_{i,j,l}(t) \equiv w_{i,j,l} \geq 0$ with $1 = \sum_{l=1}^{m} w_{i,j,l}$ is the usual backward-looking expectation, which gives $q_{i,j}^{\exp}(t) = \sum_{l=1}^{m} w_{i,j,l} q_j(k-l)$ for $t \in [k, k+1)$.

2) $q_{i,j}^{\exp}(t) = a_{i,j}(t)\int_{-T}^{-r} h_{i,j}(s)q_j(t+s)ds + (1-a_{i,j}(t))q_j^*$, where $0 < r < T$, $a_{i,j}(t) \in [0,1]$ for all $t \geq 0$,

$h_{i,j}:[-T,-r] \to \Re$ is a Lebesgue integrable function with $h_{i,j}(s) \geq 0$ for almost all $s \in [-T,-r]$ and $1 = \int_{-T}^{-r} h_{i,j}(s)ds$.

Of course, in this case it is required that $q_j(t)$ must be Lebesgue integrable and essentially bounded.

We notice the following important fact for consistent backward-looking expectations:

**FACT I:** $q_{i,j}^{\exp}(t)$ (where $j \neq i$, $i,j = 1,...,n$) is a Consistent Backward-looking expectation with respect to the Nash equilibrium point $q^* \in S$ if and only if there exist constants $0 < r \leq T$ and a function $d_{i,j}:\Re^+ \to [-1,1]$ such that:

$$q_{i,j}^{\exp}(t) = \min\left\{Q_j, \max\left\{0, q_j^* + d_{i,j}(t)\|q_j - q_j^*\|_{[t-T,t-r]}\right\}\right\}, \quad \forall t \geq 0 \tag{2.9}$$

**Proof of Fact I:** Assume first that $q_{i,j}^{\exp}(t)$ (where $j \neq i$, $i,j = 1,...,n$) is a Consistent Backward-looking expectation with respect to the Nash equilibrium point $q^* \in S$, i.e., that (2.8) holds. We distinguish the following cases.

Case 1: If $q_{i,j}^{\exp}(t) = 0$, then (2.8) implies that $q_j^* \leq \|q_j - q_j^*\|_{[t-T,t-r]}$. In this case we define $d_{i,j}(t) = -1$ and equality (2.9) holds.

Case 2: If $q_{i,j}^{\exp}(t) = Q_j$, then (2.8) implies that $Q_j - q_j^* \leq \|q_j - q_j^*\|_{[t-T,t-r]}$. In this case we define $d_{i,j}(t) = 1$ and equality (2.9) holds.

Case 3: If $q_{i,j}^{\exp}(t) \in (0, Q_j)$ and $\|q_j - q_j^*\|_{[t-T,t-r]} > 0$ then equality (2.9) holds with $d_{i,j}(t) = \text{sgn}\left(q_{i,j}^{\exp}(t) - q_j^*\right)\dfrac{\left|q_{i,j}^{\exp}(t) - q_j^*\right|}{\|q_j - q_j^*\|_{[t-T,t-r]}}$. Inequality (2.8) implies that $|d_{i,j}(t)| \leq 1$.

Case 4: If $q_{i,j}^{\exp}(t) \in (0, Q_j)$ and $\|q_j - q_j^*\|_{[t-T,t-r]} = 0$ then inequality (2.8) implies that $q_{i,j}^{\exp}(t) = q_j^*$. In this case equality (2.9) holds for arbitrary $d_{i,j}(t) \in [-1,1]$.

On the other hand, if (2.9) holds then $q_{i,j}^{\exp}(t) \in [0, Q_j]$ for all $t \geq 0$. Moreover, the reader can verify that inequality (2.8) holds. The proof is complete. ◁

For the dynamical system (2.7) we make the following assumption:

**(H):** All expectation rules $q_{i,j}^{\exp}(t)$ ($j \neq i$, $i,j = 1,...,n$) are Consistent Backward-looking expectations with respect to the Nash equilibrium point $q^* \in S$.

The previous fact shows that hypothesis (H) is equivalent to the existence of constants $0 < r \leq T$ and functions $d_{i,j}:\Re^+ \to [-1,1]$ ($j \neq i$, $i,j = 1,...,n$) such that the following equalities hold for all $i = 1,...,n$:

$$q_i(t) = \theta_i(t) \min\{Q_i, \max\{0, q_i(t - \tau_i(t))\}\}$$
$$+ (1 - \theta_i(t)) \min\left\{Q_i, \max\left\{0, \frac{ab - c_i}{2b + K_i} - \frac{b}{2b + K_i} \sum_{j \neq i} \min\left\{Q_j, \max\left\{0, q_j^* + d_{i,j}(t) \|q_j - q_j^*\|_{[t-T, t-r]}\right\}\right\}\right\}\right\} \quad (2.10)$$

In general the functions $\theta_i : \Re^+ \to [0, \Theta]$, $\tau_i : \Re^+ \to [r, T]$, $d_{i,j} : \Re^+ \to [-1, 1]$ ($j \neq i$, $i, j = 1, ..., n$) as well as the constants $0 \leq \Theta < 1$, $0 < r \leq T$ are unknown. Therefore, the dynamical system (2.10) is an uncertain dynamical system described by Functional Difference Equations (FDEs) (see [9,15-19,32]). In order to study the behavior of the solutions of (2.10) we define the dimensionless deviation variables $x_i(t) = \frac{q_i(t) - q_i^*}{Q_i}$ ($i = 1, ..., n$) and we obtain from (2.10) for $i = 1, ..., n$:

$$x_i(t) = \theta_i(t) \min\{1 - L_i, \max\{-L_i, x_i(t - \tau_i(t))\}\}$$
$$+ (1 - \theta_i(t)) \min\left\{1 - L_i, \max\left\{-L_i, M_i - L_i - R_i \sum_{j \neq i} g_{i,j} \min\left\{1, \max\left\{0, L_j + d_{i,j}(t) \|x_j\|_{[t-T, t-r]}\right\}\right\}\right\}\right\} \quad (2.11)$$

where $L_i = \frac{q_i^*}{Q_i} \in [0, 1]$, $M_i = \frac{ab - c_i}{(2b + K_i)Q_i}$, $R_i = \frac{b}{2b + K_i} > 0$, $g_{i,j} = \frac{Q_j}{Q_i} > 0$ for $j \neq i, i = 1, ..., n$ are constants which satisfy $L_i = \min\left\{1, \max\left\{0, M_i - R_i \sum_{j \neq i} g_{i,j} L_j\right\}\right\}$ for all $i = 1, ..., n$.

**Remarks and Examples about systems (2.10), (2.11):**

**a)** The reader should notice that system (2.11) is an infinite-dimensional dynamical system with state space $X$ being the normed linear space of bounded functions $x : [-T, 0] \to \Re^n$ with norm $\|x\|_X = \sup_{-T \leq \tau \leq 0} |x(\tau)|$. Indeed, by using the method of steps, given an initial condition $x_0 \in X$ and functions $\theta_i : \Re^+ \to [0, \Theta]$, $\tau_i : \Re^+ \to [r, T]$, $d_{i,j} : \Re^+ \to [-1, 1]$ ($j \neq i$, $i, j = 1, ..., n$) then one can in principle determine from (2.11) the solution $x(t) = (x_1(t), ..., x_n(t))' \in \Re^n$ for $t \in (0, r]$ with $x(\tau) = (x_1(\tau), ..., x_n(\tau))' = x_0(\tau)$ for all $\tau \in [-T, 0]$. Then we can determine from (2.11) the solution $x(t) = (x_1(t), ..., x_n(t))' \in \Re^n$ for $t \in (r, 2r]$. Continuing this way, we can determine from (2.11) the solution $x(t) = (x_1(t), ..., x_n(t))' \in \Re^n$ for $t \in (kr, (k+1)r]$, where $k$ is a positive integer. The solution is indeed bounded and exists for all $t \geq 0$, since (2.11) guarantees that $x_i(t) \in [-L_i, 1 - L_i]$ for all $t \geq 0$, $i = 1, ..., n$. The state of system (2.11) will be denoted by $x_t = \{x(\tau) : t - T \leq \tau \leq t\} \in X$ as usually in systems with delays (see [9]) and the components of the state by $x_{i,t} = \{x_i(\tau) : t - T \leq \tau \leq t\}$ for $i = 1, ..., n$.

**b)** The reader should also notice that $0 \in X$ is an equilibrium point for system (2.11). Indeed, for every functions $\theta_i : \Re^+ \to [0, \Theta]$, $\tau_i : \Re^+ \to [r, T]$, $d_{i,j} : \Re^+ \to [-1, 1]$ ($j \neq i$, $i, j = 1, ..., n$), $x_{t_0} = 0 \in X$ implies $x_t = 0 \in X$ for all $t \geq t_0$. This equilibrium point corresponds to the Nash equilibrium point $q^* \in S$ (the deviation variables have been defined by $x_i(t) = \frac{q_i(t) - q_i^*}{Q_i}$ for $i = 1, ..., n$).

**c)** All discrete-time models of the form:

$$q_i(k+1) = \theta_i(k) q_i(k) + (1 - \theta_i(k)) \min\left\{Q_i, \max\left\{0, \frac{ab - c_i}{2b + K_i} - \frac{b}{2b + K_i} \sum_{j \neq i} q_{i,j}^{\exp}(k+1)\right\}\right\}, \quad i = 1, ..., n \quad (2.12)$$

with

$$q_{i,j}^{\exp}(k+1) = a_{i,j}(k)\sum_{l=0}^{m} w_{i,j,l}(k)q_j(k-l) + (1-a_{i,j}(k))q_j^* \qquad (2.13)$$

where $k, m$ are non-negative integers, $a_{i,j}(k) \in [0,1]$ ($i, j = 1,...,n$), $\theta_i(k) \in [0, \Theta]$ ($i = 1,...,n$) with $\Theta \in [0,1)$, $w_{i,j,l}(k) \geq 0$ with $1 = \sum_{l=0}^{m} w_{i,j,l}(k)$ for all $k \geq 0$ and $l = 0,...,m$ ($i, j = 1,...,n$), are included in the uncertain model (2.10) and its equivalent expression (2.11) in the sense that for every model of the form (2.12), (2.13) one can give functions $\theta_i : \Re^+ \to [0, \Theta]$, $\tau_i : \Re^+ \to [r, T]$, $d_{i,j} : \Re^+ \to [-1,1]$ ($j \neq i$, $i, j = 1,...,n$) such that the solution of (2.10) coincides with the solution obtained by the discrete-time model (2.12), (2.13).

**d)** All continuous-time models of the form:

$$\dot{q}_i(t) = \mu_i \min\left\{ Q_i, \max\left\{ 0, \frac{ab - c_i}{2b + K_i} - \frac{b}{2b + K_i} \sum_{j \neq i} q_{i,j}^{\exp}(t) \right\} \right\} - \mu_i q_i(t), \; i = 1,...,n \qquad (2.14)$$

where $\mu_i > 0$ are constants and $q_{i,j}^{\exp}(t)$ ($j \neq i, i, j = 1,...,n$) are Consistent Backward-looking expectations with respect to the Nash equilibrium point $q^* \in S$, are included in the uncertain model (2.10). Indeed, for $t \geq r > 0$ the solution of (2.14) implies the following integral equations:

$$q_i(t) = \exp(-\mu_i r) q_i(t-r) + \int_{t-r}^{t} \exp(-\mu_i(t-\tau)) \min\left\{ Q_i, \max\left\{ 0, \frac{ab - c_i}{2b + K_i} - \frac{b}{2b + K_i} \sum_{j \neq i} q_{i,j}^{\exp}(\tau) \right\} \right\} d\tau, \; i = 1,...,n$$

From the above expression under the assumption that the mappings $t \to q_{i,j}^{\exp}(t)$ ($j \neq i, i, j = 1,...,n$) are continuous, we can conclude that for all $t \geq r$ and $i = 1,...,n$, there exists $g_i(t) \in [t-r, t]$, $i = 1,...,n$ such that

$$q_i(t) = \exp(-\mu_i r) q_i(t-r) + (1 - \exp(-\mu_i r)) \min\left\{ Q_i, \max\left\{ 0, \frac{ab - c_i}{2b + K_i} - \frac{b}{2b + K_i} \sum_{j \neq i} q_{i,j}^{\exp}(g_i(t)) \right\} \right\}, \; i = 1,...,n$$

The reader may verify that for Consistent Backward-looking expectations with respect to the Nash equilibrium point $q^* \in S$, the above model can be described by the uncertain model (2.10) with $\theta_i(t) \equiv \exp(-\mu_i r)$, $\tau_i(t) \equiv r$, $i = 1,...,n$ and $\Theta := \max_{i=1,...,n} \exp(-\mu_i r) < 1$.

The crucial question that can be posed is the question of robust asymptotic stability of the Nash equilibrium $q^* \in S$ for system (2.10) or equivalently the question of robust asymptotic stability $0 \in X$ for system (2.11). The reader can obtain rigorous definitions for robust global asymptotic stability in [15-19]. The following theorem is the main result of the present section and shows that for certain values of the parameters involved, the Nash equilibrium $q^* \in S$ is robustly globally asymptotically stable for system (2.10) in the sense that for every initial condition and for every set of functions $\theta_i : \Re^+ \to [0, \Theta]$, $\tau_i : \Re^+ \to [r, T]$, $d_{i,j} : \Re^+ \to [-1,1]$ ($j \neq i$, $i, j = 1,...,n$), the solution $q(t) = (q_1(t),..., q_n(t))' \in \Re^n$ of system (2.10) satisfies $\lim_{t \to +\infty} q(t) = q^*$.

**Theorem 2.2:** $0 \in X$ is Robustly Globally Asymptotically Stable for system (2.11), if the following set of conditions holds for each $p = 2,...,n$:

$$R_{i_1} \ldots R_{i_p}(n-1)^p < 1 \qquad (2.15)$$

for all $i_j \in \{1,...,n\}$, $i_j \neq i_k$ if $j \neq k$.

In other words, if conditions (2.15) hold then the Nash equilibrium point $q^* \in S$ is robustly Globally Asymptotically Stable with respect to all possible Consistent Backward-looking expectation rules with respect to the Nash equilibrium point $q^* \in S$, $q_{i,j}^{\exp}(t)$, $i, j = 1,...,n$, $i \neq j$. It should be noticed that conditions (2.15) are more demanding inequalities than other stability conditions in the literature. However, this is expected since conditions (2.15) are sufficient conditions for global asymptotic stability for the uncertain model (2.10) which contains many models studied in the literature as special cases.

Conditions (2.15) are termed as small-gain conditions in Mathematical Control Theory (see [17,18,19]). For $n = 2$ conditions (2.15) are equivalent to the inequality:

$$R_1 R_2 < 1$$

For $n = 3$, conditions (2.15) are equivalent to the following four inequalities:

$$4R_1 R_2 < 1, \; 4R_1 R_3 < 1, \; 4R_2 R_3 < 1, \; 8R_1 R_2 R_3 < 1 \tag{2.16}$$

For $n = 4$, conditions (2.15) are equivalent to the following eleven inequalities:

$$9R_1 R_2 < 1, \; 9R_1 R_3 < 1, \; 9R_1 R_4 < 1, \; 9R_2 R_3 < 1, \; 9R_2 R_4 < 1, \; 9R_3 R_4 < 1$$

$$27 R_1 R_2 R_3 < 1, \; 27 R_1 R_2 R_4 < 1, \; 27 R_1 R_3 R_4 < 1, \; 27 R_2 R_3 R_4 < 1$$

$$81 R_1 R_2 R_3 R_4 < 1$$

The proof of Theorem 2.2 relies heavily on recent results on dynamical systems (see [19]) and techniques developed for delay systems (see [18,19]) and is provided at the Appendix. An interesting corollary for the Cournot oligopoly game is given next.

**Corollary 2.3:** *If conditions (2.15) hold for each $p = 2,...,n$ and for all $i_j \in \{1,...,n\}$, $i_j \neq i_k$ if $j \neq k$ then the Nash equilibrium point $q^* \in S$ is unique for the game described by (2.3), (2.4) and (2.5).*

The reader should notice that Brouwer's fixed point theorem guarantees the existence of the Nash equilibrium $q^* \in S$ but does not guarantee uniqueness.

**Proof of Corollary 2.3:** Suppose that there exists $q^{**} \in S$ with $q_i^{**} = \min\left\{Q_i, \max\left\{0, \frac{ab - c_i}{2b + K_i} - \frac{b}{2b + K_i} \sum_{j \neq i} q_j^{**}\right\}\right\}$ for $i = 1,...,n$ and $q^{**} \neq q^*$. This implies that $y_i := \frac{q_i^{**} - q_i^*}{Q_i} = \min\left\{1 - L_i, \max\left\{-L_i, M_i - L_i - R_i \sum_{j \neq i} g_{i,j}(L_j + y_j)\right\}\right\}$ for $i = 1,...,n$. Using the previous equalities, the reader can verify that the solution of (2.11) with initial condition $x_0 = Py$, where $Py := \{y : -T \leq \tau \leq 0\} \in X$, $y = (y_1,..., y_n)'$, corresponding to the constant inputs

$$\theta_i(t) \equiv 0, \; d_{i,j}(t) := \operatorname{sgn}(y_j) = \begin{cases} 1 & \text{if } y_j > 0 \\ 0 & \text{if } y_j = 0 \\ -1 & \text{if } y_j < 0 \end{cases}, \; j \neq i, \; i, j = 1,...,n$$

satisfies $x_t = Py$ for all $t \geq 0$ ($\tau_i : \mathfrak{R}^+ \to [r, T]$, $i = 1,...,n$ are irrelevant) and consequently we cannot have $\lim_{t \to +\infty} x(t) = 0$. This is impossible according to Theorem 2.2. The proof is complete. ◁

# 3. Extension to the General Case of Dynamic Games

The approach described in the previous section for the Cournot oligopoly game can be extended to any strategic game.

Consider a strategic game with $n$ players and $S_i \subseteq \Re^{k_i}$ ($i = 1,...,n$) being the action space for each one of the players. We assume that the best reply mapping for each one of the players is a function $f_i : S_1 \times ... \times S_{i-1} \times S_{i+1} \times ... \times S_n \to S_i$ for $1 < i < n$, $n \geq 3$ and $f_1 : S_2 \times ... \times S_n \to S_1$, $f_n : S_1 \times ... \times S_{n-1} \to S_n$, satisfying the following inequalities:

$$\pi_i(q_i, q_{-i}) < \pi_i(f_i(q_{-i}), q_{-i}), \text{ for all } q_i \in S_i \text{ with } q_i \neq f_i(q_{-i}), i = 1,...,n \quad (3.1)$$

where $\pi_i(q_i, q_{-i})$ is the payoff function of the i-*th* player.

We assume the existence of a Nash equilibrium $q^* \in S$ for the game, where $S := S_1 \times ... \times S_n$ is the outcome space for the game, i.e., there exists $q^* = (q_1^*,...,q_n^*) \in S$ such that

$$q_i^* = f_i(q_{-i}^*), i = 1,...,n \quad (3.2)$$

The existence of a Nash equilibrium can be guaranteed by Brouwer's fixed point theorem when all action spaces $S_i \subseteq \Re^{k_i}$ ($i = 1,...,n$) are compact and convex and when all the best reply mappings $f_i : S_1 \times ... \times S_{i-1} \times S_{i+1} \times ... \times S_n \to S_i$ for $1 < i < n$, $n \geq 3$ and $f_1 : S_2 \times ... \times S_n \to S_1$, $f_n : S_1 \times ... \times S_{n-1} \to S_n$ are continuous mappings.

Next we assume that $S_i \subseteq \Re^{k_i}$ ($i = 1,...,n$) are closed convex sets and that the dynamics of the game are described in continuous time as follows:

- every player forms an expectation for the behavior of all other players at each time $t \geq 0$: the expectation of the i-*th* player for the production level of the j-*th* player at time $t \geq 0$ will be denoted by $q_{i,j}^{\exp}(t) \in S_j$ ($j \neq i$, $i, j = 1,...,n$),

- every player determines her action as a convex combination of a past action and the best reply response based on the expectations for the behavior of all other players at each time $t \geq 0$, i.e.,

$$q_1(t) = \theta_1(t) \Pr_{S_1}(q_1(t - \tau_1(t))) + (1 - \theta_1(t)) f_1(q_{1,2}^{\exp}(t),...,q_{1,n}^{\exp}(t))$$
$$\vdots$$
$$q_i(t) = \theta_i(t) \Pr_{S_i}(q_i(t - \tau_i(t))) + (1 - \theta_i(t)) f_i(q_{i,1}^{\exp}(t),...q_{i,i-1}^{\exp}(t), q_{i,i+1}^{\exp}(t),...,q_{i,n}^{\exp}(t)) \quad (3.3)$$
$$\vdots$$
$$q_n(t) = \theta_n(t) \Pr_{S_n}(q_n(t - \tau_n(t))) + (1 - \theta_n(t)) f_n(q_{n,1}^{\exp}(t),...,q_{n,n-1}^{\exp}(t))$$

where $\theta_i : \Re^+ \to [0, \Theta]$, $\tau_i : \Re^+ \to [r, T]$, $i = 1,...,n$ are in general unknown functions, $0 \leq \Theta < 1$, $0 < r \leq T$ are constants (in general unknown),

- all expectation rules $q_{i,j}^{\exp}(t) \in S_j$ ($j \neq i, i = 1,...,n$), are Consistent Backward-looking expectations with respect to the Nash equilibrium point $q^* \in S$, i.e., there exist constants $0 < r \leq T$ such that:

$$\left| q_{i,j}^{\exp}(t) - q_j^* \right| \leq \sup_{t-T \leq \tau \leq t-r} \left| q_j(\tau) - q_j^* \right| = \left\| q_j - q_j^* \right\|_{[t-T, t-r]}, \text{ for all } t \geq 0 \quad (3.4)$$

We notice the following fact for consistent backward-looking expectations:

**FACT II:** *Suppose that $S_j \subseteq \Re^{k_j}$ is a closed convex set. $q_{i,j}^{\exp}(t)$ ($j \neq i$, $i,j = 1,...,n$) is a Consistent Backward-looking expectation with respect to the Nash equilibrium point $q^* \in S$ if and only if there exist constants $0 < r \leq T$ and a function $d_{i,j} : \Re^+ \to \{d \in \Re^{k_j} : |d| \leq 1\}$ such that:*

$$q_{i,j}^{\exp}(t) = \Pr\nolimits_{S_j}\left(q_j^* + d_{i,j}(t)\|q_j - q_j^*\|_{[t-T,t-r]}\right), \quad \forall t \geq 0 \tag{3.5}$$

Proof of Fact II: Indeed, using the fact that $|\Pr_U(x) - \Pr_U(y)| \leq |x - y|$ for every $x, y \in \Re^n$, where $U \subseteq \Re^n$ is a closed convex set, one can verify that for every $d_{i,j} : \Re^+ \to \{d \in \Re^{k_j} : |d| \leq 1\}$ the function $q_{i,j}^{\exp}(t)$ defined by (3.5) satisfies (3.4) and $q_{i,j}^{\exp}(t) \in S_j$ for all $t \geq 0$. Hence it is a Consistent Backward-looking expectation with respect to the Nash equilibrium point $q^* \in S$. On the other hand, if $q_{i,j}^{\exp}(t) \in S_j$ is a Consistent Backward-looking expectation with respect to the Nash equilibrium point $q^* \in S$ satisfying (3.4) for all $t \geq 0$ then the function defined by:

$$d_{i,j}(t) = \frac{1}{\|q_j - q_j^*\|_{[t-T,t-r]}} \left(q_{i,j}^{\exp}(t) - q_j^*\right), \text{ if } \|q_j - q_j^*\|_{[t-T,t-r]} > 0$$

$$d_{i,j}(t) = 0, \text{ if } \|q_j - q_j^*\|_{[t-T,t-r]} = 0$$

satisfies $d_{i,j}(t) \in \{d \in \Re^{k_j} : |d| \leq 1\}$. Moreover, (3.5) holds for all $t \geq 0$. The proof is complete. ◁

Fact II shows that if all action spaces $S_j \subseteq \Re^{k_j}$ ($j = 1,...,n$) are closed convex sets then there exist constants $0 < r \leq T$, $0 \leq \Theta < 1$ and functions $\theta_i : \Re^+ \to [0, \Theta]$, $\tau_i : \Re^+ \to [r, T]$, $d_{i,j} : \Re^+ \to \{d \in \Re^{k_j} : |d| \leq 1\}$ ($j \neq i$, $i, j = 1,...,n$) such that:

$$q_1(t) = \theta_1(t)\Pr\nolimits_{S_1}(q_1(t - \tau_1(t))) + (1 - \theta_1(t))f_1\left(\Pr\nolimits_{S_2}\left(q_2^* + d_{1,2}(t)\|q_2 - q_2^*\|_{[t-T,t-r]}\right),...,\Pr\nolimits_{S_n}\left(q_n^* + d_{1,n}(t)\|q_n - q_n^*\|_{[t-T,t-r]}\right)\right),$$
$$\vdots$$
$$q_n(t) = \theta_n(t)\Pr\nolimits_{S_n}(q_n(t - \tau_n(t))) + (1 - \theta_n(t))f_n\left(\Pr\nolimits_{S_1}\left(q_1^* + d_{n,1}(t)\|q_1 - q_1^*\|_{[t-T,t-r]}\right),...,\Pr\nolimits_{S_{n-1}}\left(q_{n-1}^* + d_{n,n-1}(t)\|q_{n-1} - q_{n-1}^*\|_{[t-T,t-r]}\right)\right)$$
$$\tag{3.6}$$

In general, the constants $0 < r \leq T$, $0 \leq \Theta < 1$ and the functions $\theta_i : \Re^+ \to [0, \Theta]$, $\tau_i : \Re^+ \to [r, T]$, $d_{i,j} : \Re^+ \to \{d \in \Re^{k_j} : |d| \leq 1\}$ ($j \neq i$, $i, j = 1,...,n$) are unknown. Therefore, the dynamical system (3.6) is an uncertain dynamical system described by Functional Difference Equations (FDEs) (see [9,15-19,32]). In order to study the behavior of the solutions of (3.6) we define the deviation variables $x_i(t) = q_i(t) - q_i^*$ ($i = 1,...,n$) and we obtain from (3.6):

$$x_1(t) = \theta_1(t)\left(\Pr\nolimits_{S_1}\left(x_1(t - \tau_1(t)) + q_1^*\right) - q_1^*\right)$$
$$+ (1 - \theta_1(t))\left(f_1\left(\Pr\nolimits_{S_2}\left(q_2^* + d_{1,2}(t)\|x_2\|_{[t-T,t-r]}\right),...,\Pr\nolimits_{S_n}\left(q_n^* + d_{1,n}(t)\|x_n\|_{[t-T,t-r]}\right)\right) - q_1^*\right),$$
$$\vdots \tag{3.7}$$
$$x_n(t) = \theta_n(t)\left(\Pr\nolimits_{S_n}\left(x_n(t - \tau_n(t)) + q_n^*\right) - q_n^*\right)$$
$$+ (1 - \theta_n(t))\left(f_n\left(\Pr\nolimits_{S_1}\left(q_1^* + d_{n,1}(t)\|x_1\|_{[t-T,t-r]}\right),...,\Pr\nolimits_{S_{n-1}}\left(q_{n-1}^* + d_{n,n-1}(t)\|x_{n-1}\|_{[t-T,t-r]}\right)\right) - q_n^*\right).$$

Finally, we assume that there exist functions $\tilde{\gamma}_{i,j} \in N$ ($j \neq i$, $i,j = 1,...,n$) such that the following inequalities hold for all $q \in S$:

$$\left| f_i(q_{-i}) - q_i^* \right| \leq \max_{j \neq i} \tilde{\gamma}_{i,j}\left( \left| q_j - q_j^* \right| \right), \; i = 1,...,n \tag{3.8}$$

Using again the fact that $\left| \Pr_U(x) - \Pr_U(y) \right| \leq \left| x - y \right|$ for every $x, y \in \mathfrak{R}^n$, where $U \subseteq \mathfrak{R}^n$ is a closed convex set and inequalities (3.8), we obtain from (3.7) for all $t \geq 0$ and $\mu > \Theta$:

$$\left| x_i(t) \right| \leq \max\left\{ \mu \left\| x_i \right\|_{[t-T,t-r]}, \max_{j \neq i} \frac{\mu - \mu\Theta}{\mu - \Theta} \tilde{\gamma}_{i,j}\left( \left\| x_j \right\|_{[t-T,t-r]} \right) \right\}, \; i = 1,...,n \tag{3.9}$$

**Remarks and Examples about systems (3.7), (3.8):**

**a)** The reader should notice that system (3.7) is an infinite-dimensional dynamical system with state space $\mathcal{X}$ being the normed linear space of bounded functions $x:[-T,0] \to \mathfrak{R}^N$, where $N = k_1 + ... + k_n$ with norm $\left\| x \right\|_\mathcal{X} = \sup_{-T \leq \tau \leq 0} \left| x(\tau) \right|$. Indeed, by using the method of steps, given an initial condition $x_0 \in \mathcal{X}$ and functions $\theta_i : \mathfrak{R}^+ \to [0,\Theta]$, $\tau_i : \mathfrak{R}^+ \to [r,T]$, $d_{i,j} : \mathfrak{R}^+ \to \left\{ d \in \mathfrak{R}^{k_j} : \left| d \right| \leq 1 \right\}$ ($j \neq i$, $i,j = 1,...,n$) then one can in principle determine from (3.8) the solution $x(t) = (x_1(t),...,x_n(t))' \in \mathfrak{R}^m$ for $t \in (0,r]$ with $x(\tau) = (x_1(\tau),...,x_n(\tau))' = x_0(\tau)$ for all $\tau \in [-T,0]$. Then we can determine from (3.8) the solution $x(t) = (x_1(t),...,x_n(t))' \in \mathfrak{R}^N$ for $t \in (r,2r]$. Continuing this way, we can determine from (3.8) the solution $x(t) = (x_1(t),...,x_n(t))' \in \mathfrak{R}^N$ for $t \in (kr,(k+1)r]$, where $k$ is a positive integer. The state of system (3.8) will be denoted by $x_t = \{x(\tau) : t-T \leq \tau \leq t\} \in \mathcal{X}$ as usually in systems with delays (see [9]) and the components of the state by $x_{i,t} = \{x_i(\tau) : t-T \leq \tau \leq t\}$ for $i = 1,...,n$. The solution exists for all $t \geq 0$ and satisfies $x_t = \{x(\tau) : t-T \leq \tau \leq t\} \in \mathcal{X}$ for all $t \geq 0$. To see this, notice that (3.9) implies the existence of a function $G \in N$ such that:

$$\sup_{t \in [0,r]} \left| x(t) \right| = \left\| x \right\|_{[0,r]} \leq G\left( \left\| x \right\|_{[-T,0]} \right) \tag{3.10}$$

Without loss of generality we may assume that $G(s) \geq s$ for all $s \geq 0$. Inequality (3.10) implies that

$$\left\| x \right\|_{[t-T,t]} \leq G\left( \left\| x \right\|_{[-T,0]} \right), \text{ for all } t \in [0,r] \text{ and } \left\| x \right\|_{[r-T,r]} \leq G\left( \left\| x \right\|_{[-T,0]} \right) \tag{3.11}$$

Working in this way and using induction we may establish that for every positive integer $k > 0$ it holds that

$$\left\| x \right\|_{[t-T,t]} \leq G^{(k)}\left( \left\| x \right\|_{[-T,0]} \right), \text{ for all } t \in [0,kr] \tag{3.12}$$

where $G^{(k)}(s) := \underbrace{G \circ ... \circ G}_{k \text{ times}}(s)$. Therefore (3.12) implies that

$$\left\| x_t \right\|_\mathcal{X} = \left\| x \right\|_{[t-T,t]} \leq G^{(1+[t/r])}\left( \left\| x \right\|_{[-T,0]} \right) = G^{(1+[t/r])}\left( \left\| x_0 \right\|_\mathcal{X} \right), \text{ for all } t \geq 0 \tag{3.13}$$

where $[t/r]$ denotes the integer part of $t/r$.

**b)** The reader should also notice that $0 \in \mathcal{X}$ is an equilibrium point for system (3.7). Indeed, for any functions $\theta_i : \mathfrak{R}^+ \to [0,\Theta]$, $\tau_i : \mathfrak{R}^+ \to [r,T]$, $d_{i,j} : \mathfrak{R}^+ \to \left\{ d \in \mathfrak{R}^{k_j} : \left| d \right| \leq 1 \right\}$ ($j \neq i$, $i,j = 1,...,n$), $x_{t_0} = 0 \in \mathcal{X}$ implies $x_t = 0 \in \mathcal{X}$ for all $t \geq t_0$. This equilibrium point corresponds to the Nash equilibrium point $q^* \in S$ (noting that the deviation variables have been defined by $x_i(t) = q_i(t) - q_i^*$ for $i = 1,...,n$).

**c)** All discrete-time models of the form:

$$q_1(k+1) = \theta_1(k)q_1(k) + (1-\theta_1(k))f_1(q_{1,2}^{\exp}(k+1),...,q_{1,n}^{\exp}(k+1))$$
$$\vdots$$
$$q_i(k+1) = \theta_i(k)q_i(k) + (1-\theta_i(k))f_i(q_{i,1}^{\exp}(k+1),...q_{i,i-1}^{\exp}(k+1), q_{i,i+1}^{\exp}(k+1),...,q_{i,n}^{\exp}(k+1)) \quad (3.14)$$
$$\vdots$$
$$q_n(k+1) = \theta_n(k)q_n(k) + (1-\theta_n(k))f_n(q_{n,1}^{\exp}(k+1),...,q_{n,n-1}^{\exp}(k+1))$$

with $q_{i,j}^{\exp}(k+1)$ given by (2.13), where $k,m$ are non-negative integers, $a_{i,j}(k) \in [0,1]$ ($i,j=1,...,n$), $\theta_i(k) \in [0,\Theta]$ ($i=1,...,n$) with $\Theta \in [0,1)$, $w_{i,j,l}(k) \geq 0$ with $1 = \sum_{l=0}^{m} w_{i,j,l}(k)$ for all $k \geq 0$ and $l=0,...,m$ ($i,j=1,...,n$), are included in the uncertain model (3.6) and its equivalent expression (3.7) in the sense that for every model of the form (3.14), (2.13), one can give functions $\theta_i : \Re^+ \to [0,\Theta]$, $\tau_i : \Re^+ \to [r,T]$, $d_{i,j} : \Re^+ \to \{d \in \Re^{k_j} : |d| \leq 1\}$ ($j \neq i$, $i,j=1,...,n$) such that the solution of (3.6) coincides with the solution obtained by the discrete-time model (3.14), (2.13).

**d)** Similarly, as shown in previous section, if $S_j \subseteq [0,+\infty)$ for all $j=1,...,n$ and if all expectation rules $q_{i,j}^{\exp}(t)$ ($j \neq i$, $i,j=1,...,n$) are Consistent Backward-looking expectations with respect to the Nash equilibrium point $q^* \in S$ and all mappings $t \to q_{i,j}^{\exp}(t)$ ($j \neq i$, $i,j=1,...,n$) are continuous, then all continuous-time models of the form:

$$\dot{q}_1(t) = \mu_1 \left( f_1(q_{1,2}^{\exp}(t),...,q_{1,n}^{\exp}(t)) - q_1(t) \right)$$
$$\vdots$$
$$\dot{q}_n(t) = \mu_n \left( f_n(q_{n,1}^{\exp}(t),...,q_{n,n-1}^{\exp}(t)) - q_n(t) \right)$$

where $\mu_i > 0$ are constants, are included in the uncertain model (3.6).

**e)** The reader should notice that no continuity assumption is made for the best reply mappings of the players $f_i : S_1 \times ... \times S_{i-1} \times S_{i+1} \times ... \times S_n \to S_i$ for $1 < i < n$, $n \geq 3$ and $f_1 : S_2 \times ... \times S_n \to S_1$, $f_n : S_1 \times ... \times S_{n-1} \to S_n$. Moreover, we have not assumed that the action spaces $S_j \subseteq \Re^{k_j}$ ($j=1,...,n$) are compact sets: we simply require that the action spaces are closed, convex sets. However, we have assumed the existence of a Nash equilibrium point $q^* \in S$ and the existence of functions $\tilde{\gamma}_{i,j} \in N$ ($j \neq i$, $i=1,...,n$) satisfying (3.8).

The crucial question that can be posed is the question of robust asymptotic stability of the Nash equilibrium $q^* \in S$ for system (3.6) or equivalently the question of robust asymptotic stability $0 \in X$ for system (3.7). The following theorem is the main result of this section and shows that robust global stability can be determined by the functions $\tilde{\gamma}_{i,j} \in N$ ($j \neq i$, $i=1,...,n$) satisfying (3.8).

**Theorem 3.1:** $0 \in X$ *is Robustly Globally Asymptotically Stable for system (3.7), if there exists $\omega > 1$ such that the following set of conditions holds for each $p = 2,...,n$:*

$$\left( \gamma_{i_1,i_2} \circ \gamma_{i_2,i_3} \circ ... \circ \gamma_{i_p,i_1} \right)(s) < s, \quad \forall s > 0 \quad (3.15)$$

*for all $i_j \in \{1,...,n\}$, $i_j \neq i_k$ if $j \neq k$, where $\gamma_{i,j}(s) := \omega \tilde{\gamma}_{i,j}(\omega s)$.*

In other words, if conditions (3.15) hold then the Nash equilibrium point $q^* \in S$ is robustly Globally Asymptotically Stable with respect to all possible Consistent Backward-looking expectation rules with respect to the Nash equilibrium point $q^* \in S$, $q_{i,j}^{\exp}(t)$, $i, j = 1,...,n$, $i \neq j$.

Conditions (3.15) are termed as cyclic small-gain conditions in Mathematical Control Theory (see [17,18,19]). For $n = 2$ conditions (3.15) are equivalent to the inequalities:

$$\gamma_{1,2}(\gamma_{2,1}(s)) < s \text{ and } \gamma_{2,1}(\gamma_{1,2}(s)) < s, \forall s > 0$$

For $n = 3$, conditions (3.15) are equivalent to the following twelve inequalities for all $s > 0$:

$$\gamma_{1,2}(\gamma_{2,1}(s)) < s, \gamma_{2,1}(\gamma_{1,2}(s)) < s$$

$$\gamma_{1,3}(\gamma_{3,1}(s)) < s, \gamma_{3,1}(\gamma_{1,3}(s)) < s$$

$$\gamma_{2,3}(\gamma_{3,2}(s)) < s, \gamma_{3,2}(\gamma_{2,3}(s)) < s$$

$$\gamma_{1,2}(\gamma_{2,3}(\gamma_{3,1}(s))) < s, \gamma_{3,1}(\gamma_{1,2}(\gamma_{2,3}(s))) < s, \gamma_{2,3}(\gamma_{3,1}(\gamma_{1,2}(s))) < s$$

$$\gamma_{2,1}(\gamma_{1,3}(\gamma_{3,2}(s))) < s, \gamma_{1,3}(\gamma_{3,2}(\gamma_{2,1}(s))) < s, \gamma_{3,2}(\gamma_{2,1}(\gamma_{1,3}(s))) < s$$

The reader should notice that many of the above inequalities are equivalent. For example, for $n = 2$ the inequality $\gamma_{1,2}(\gamma_{2,1}(s)) < s$, for all $s > 0$ implies the inequality $\gamma_{2,1}(\gamma_{1,2}(s)) < s$ for all $s > 0$. Similarly, for the case $n = 3$ the following five inequalities $\gamma_{1,2}(\gamma_{2,1}(s)) < s$, $\gamma_{1,3}(\gamma_{3,1}(s)) < s$, $\gamma_{2,3}(\gamma_{3,2}(s)) < s$, $\gamma_{1,2}(\gamma_{2,3}(\gamma_{3,1}(s))) < s$, $\gamma_{2,1}(\gamma_{1,3}(\gamma_{3,2}(s))) < s$ for all $s > 0$, imply all twelve inequalities which express conditions (3.15) in this case.

It should be noticed that for the Cournot oligopoly game studied in the previous section, the best reply mappings $f_i$ ($i = 1,...,n$) are defined by (2.3). Consequently, using the convexity of the sets $S_i = [0, Q_i]$ ($i = 1,...,n$), we obtain the following inequalities for $i = 1,...,n$:

$$\left| f_i(q_{-i}) - q_i^* \right| \leq \frac{b}{2b + K_i} \sum_{j \neq i} \left| q_j - q_j^* \right| \leq \frac{b}{2b + K_i} (n-1) \max_{j \neq i} \left| q_j - q_j^* \right|$$

The above inequalities imply that inequalities (3.8) hold with $\tilde{\gamma}_{i,j}(s) := R_i(n-1)s$, where $R_i := \frac{b}{2b + K_i}$. Theorem 3.1 and the above definitions guarantee robust global asymptotic stability of the Nash equilibrium provided that there exists $\omega > 1$ such that the following set of conditions holds for each $p = 2,...,n$:

$$R_{i_1} R_{i_2} ... R_{i_p} (n-1)^p \omega^{2p} < 1$$

for all $i_j \in \{1,...,n\}$, $i_j \neq i_k$ if $j \neq k$. Conditions (2.15) are necessary and sufficient conditions for the existence of a (sufficiently small) constant $\omega > 1$ satisfying the above inequalities for each $p = 2,...,n$ and for all $i_j \in \{1,...,n\}$, $i_j \neq i_k$ if $j \neq k$. Therefore, we conclude that Theorem 2.2 is a special case of Theorem 3.1. A more careful analysis similar to the above analysis reveals that the Nash equilibrium for the Cournot oligopoly game described in Section 2 will be asymptotically stable provided that there exist $n$ sets of positive real numbers $A_i = \{a_{i,j}, j \neq i\}$ ($i = 1,...,n$) with $\sum_{j \neq i} x_j \leq \max_{j \neq i}(a_{i,j} x_j)$ for all $x = (x_1,...,x_n)' \in (\Re^+)^n$ and $i = 1,...,n$ such that the following set of conditions holds for each $p = 2,...,n$:

$$a_{i_1,i_2} a_{i_2,i_3} ... a_{i_p,i_1} R_{i_1} R_{i_2} ... R_{i_p} < 1$$

for all $i_j \in \{1,...,n\}$, $i_j \neq i_k$ if $j \neq k$. The above conditions are less restrictive than conditions (2.15); indeed, conditions (2.15) are implied by the above conditions for the special case $a_{i,j} = n-1$ for all $i, j = 1,...,n$ with $j \neq i$. For example, for $n = 3$, the above small-gain conditions are equivalent to the existence of $\varepsilon_1, \varepsilon_2, \varepsilon_3 > 0$ such that:

$$R_1 R_2 (1+\varepsilon_1)(1+\varepsilon_2) < 1$$
$$R_1 R_3 (1+\varepsilon_1^{-1})(1+\varepsilon_3) < 1$$
$$R_2 R_3 (1+\varepsilon_2^{-1})(1+\varepsilon_3^{-1}) < 1$$
$$R_1 R_2 R_3 (1+\varepsilon_1)(1+\varepsilon_2^{-1})(1+\varepsilon_3) < 1$$
$$R_1 R_2 R_3 (1+\varepsilon_1^{-1})(1+\varepsilon_2)(1+\varepsilon_3^{-1}) < 1$$

For the above inequalities we have used $a_{1,2} = 1+\varepsilon_1$, $a_{1,3} = 1+\varepsilon_1^{-1}$, $a_{2,1} = 1+\varepsilon_2$, $a_{2,3} = 1+\varepsilon_2^{-1}$, $a_{3,1} = 1+\varepsilon_3$ and $a_{3,2} = 1+\varepsilon_3^{-1}$. By selecting $\varepsilon_1 = \varepsilon_2 = \varepsilon_3 = 1$, we obtain inequalities (2.16).

It should be emphasized that the parameters $T \geq r > 0$ which are involved in the definition of the Consistent Backward-looking expectation (Definition 2.1), play no role in the small-gain conditions. Consequently, the small-gain conditions can help us to decide whether the Nash equilibrium point is robustly stable **without any knowledge of the expectation rules**. The small-gain conditions (3.15) demand knowledge of the Nash equilibrium point $q^* \in S$ and the best reply mappings $f_i : S_1 \times ... \times S_{i-1} \times S_{i+1} \times ... \times S_n \to S_i$ for $1 < i < n$, $n \geq 3$ and $f_1 : S_2 \times ... \times S_n \to S_1$, $f_n : S_1 \times ... \times S_{n-1} \to S_n$ for which inequalities (3.8) hold.

The proof of Theorem 3.1 relies heavily on recent results on dynamical systems (see [19]) and techniques developed for time-delay systems (see [18,19]) and is provided in the Appendix. An interesting corollary is given next.

**Corollary 3.2:** *If there exists $\omega > 1$ such that conditions (3.15) hold for each $p = 2,...,n$ and for all $i_j \in \{1,...,n\}$, $i_j \neq i_k$ if $j \neq k$ then the Nash equilibrium point $q^* \in S$ is unique.*

The proof of Corollary 3.2 is exactly the same with the proof of Corollary 2.3: we show that the existence of an additional Nash equilibrium $q^{**} \in S$ contradicts robust global asymptotic stability of $0 \in X$ for system (3.7). Details are left to the reader.

Using Corollary 3.2 we may obtain conditions for uniqueness for a fixed point. Indeed, we have:

**Corollary 3.3:** *Let $S_i \subseteq \Re^{k_i}$ ($i = 1,...,n$) be closed, convex sets and let functions $f_i : S_1 \times ... \times S_{i-1} \times S_{i+1} \times ... \times S_n \to S_i$ for $1 < i < n$, $n \geq 3$ and $f_1 : S_2 \times ... \times S_n \to S_1$, $f_n : S_1 \times ... \times S_{n-1} \to S_n$ for which there exists $q^* = (q_1^*,...,q_n^*) \in S$, where $S := S_1 \times ... \times S_n$ satisfying (3.3). Furthermore, suppose that there exist functions $\tilde{\gamma}_{i,j} \in N$ ($j \neq i$, $i, j = 1,...,n$) such that inequalities (3.8) hold for all $q \in S$ and that there exists $\omega > 1$ such that conditions (3.15) hold for each $p = 2,...,n$ and for all $i_j \in \{1,...,n\}$, $i_j \neq i_k$ if $j \neq k$. Then $q^* = (q_1^*,...,q_n^*) \in S$ is the unique fixed point of the mapping $S \ni q \to F(q) := (f_1(q_{-1}),...,f_n(q_{-n})) \in S$.*

The reader should notice that Corollary 3.3 does not guarantee the existence of a fixed point for the mapping $S \ni q \to F(q) := (f_1(q_{-1}),...,f_n(q_{-n})) \in S$. Corollary 3.3 can be used in conjunction with classical fixed-point theorems (e.g., Brouwer's fixed point theorem when all action spaces $S_i \subseteq \Re^{k_i}$ ($i = 1,...,n$) are compact and convex and when the mapping $S \ni q \to F(q) := (f_1(q_{-1}),...,f_n(q_{-n})) \in S$ is continuous) in order to guarantee uniqueness of the fixed point.

# 4. Remarks on the Rational Expectation Case

It is clear from Definition 2.1 that rational expectations $q_{i,j}^{\exp}(t) = q_j(t)$ are not necessarily consistent backward-looking expectations with respect to the Nash equilibrium point $q^* \in S$. Therefore, the case of rational expectations is not necessarily covered by the results of the previous sections. This point motivates the following definition for the strategic game described in the previous section.

**Definition 4.1:** *An expectation rule $q_{i,j}^{\exp}(t)$ ($j \neq i$, $i, j = 1,...,n$) is called a Rational-Consistent Backward-looking expectation with respect to the Nash equilibrium point $q^* \in S$ if there exists a constant $0 < T$ such that:*

$$\left| q_{i,j}^{\exp}(t) - q_j^* \right| \leq \sup_{t-T \leq \tau \leq t} \left| q_j(\tau) - q_j^* \right| = \left\| q_j - q_j^* \right\|_{[t-T,t]}, \text{ for all } t \geq 0 \quad (4.1)$$

Clearly, rational expectations are Rational-Consistent Backward-looking expectations with respect to the Nash equilibrium point $q^* \in S$. Moreover, a Consistent Backward-looking expectation (in the sense of Definition 2.1) is a Rational-Consistent Backward-looking expectation with respect to the Nash equilibrium point $q^* \in S$.

Can we consider system (3.3) where all expectation rules $q_{i,j}^{\exp}(t) \in S_j$ ($j \neq i$, $i, j = 1,...,n$), are Rational-Consistent Backward-looking expectations with respect to the Nash equilibrium point $q^* \in S$? The key mathematical problem that arises in this case is whether we can obtain a well-defined dynamical system: Remark (a) in previous section does not apply. However, we can extend the analysis of the previous section under the following hypothesis:

**(H')** *There exist $m$ index sets $J_l \subseteq \{1,...,n\}$, $l = 1,...,m$ with $J_l \cap J_k = \emptyset$ for $l \neq k$ and $\underset{l=1,...,m}{\cup} J_l = \{1,...,n\}$ such that:*

*All players with $i \in J_m$ are using Consistent Backward-looking expectations with respect to the Nash equilibrium point $q^* \in S$.*

*Moreover, for every $k = 1,...,m-1$, the following statement holds:*

*All players with $i \in J_k$ are using Rational-Consistent Backward-looking expectations, $q_{i,j}^{\exp}(t) \in S_j$ if $j \in J_l$ with $l > k$ and Consistent Backward-looking expectations, $q_{i,j}^{\exp}(t) \in S_j$ if otherwise.*

Indeed, if hypothesis (H') holds then system (3.3) is expressed in deviation variables $x_i(t) = q_i(t) - q_i^*$ ($i = 1,...,n$) by the following equations:

$$\begin{aligned}
x_1(t) &= \theta_1(t)\left(\Pr_{S_1}\left(x_1(t - \tau_1(t)) + q_1^*\right) - q_1^*\right) \\
&+ (1 - \theta_1(t))\left(f_1\left(\Pr_{S_2}\left(q_2^* + d_{1,2}(t)s_{1,2}(t)\right),...,\Pr_{S_n}\left(q_n^* + d_{1,n}(t)s_{1,n}(t)\right)\right) - q_1^*\right), \\
&\vdots \\
x_n(t) &= \theta_n(t)\left(\Pr_{S_n}\left(x_n(t - \tau_n(t)) + q_n^*\right) - q_n^*\right) \\
&+ (1 - \theta_n(t))\left(f_n\left(\Pr_{S_1}\left(q_1^* + d_{n,1}(t)s_{n,1}(t)\right),...,\Pr_{S_{n-1}}\left(q_{n-1}^* + d_{n,n-1}(t)s_{n,n-1}(t)\right)\right) - q_n^*\right), \\
s_{i,j}(t) &:= \left\| x_j \right\|_{[t-T,t]} \quad \text{if} \quad i \in J_k \text{ and } j \in J_l \text{ with } l > k \\
s_{i,j}(t) &:= \left\| x_2 \right\|_{[t-T,t-r]} \quad \text{if} \quad i \in J_k \text{ and } j \in J_l \text{ with } l \leq k
\end{aligned} \quad (4.2)$$

where $\theta_i : \Re^+ \to [0, \Theta]$, $\tau_i : \Re^+ \to [r, T]$, $d_{i,j} : \Re^+ \to \left\{ d \in \Re^{k_j} : |d| \leq 1 \right\}$ ($j \neq i$, $i, j = 1,...,n$).

Let us explain next why system (4.2) is an infinite-dimensional dynamical system with state space $X$ being the normed linear space of bounded functions $x:[-T,0] \to \Re^N$, where $N = k_1 + ... + k_n$ with norm $\|x\|_X = \sup_{-T \leq \tau \leq 0} |x(\tau)|$. Indeed, by using the method of steps, given an initial condition $x_0 \in X$ and functions $\theta_i : \Re^+ \to [0, \Theta]$, $\tau_i : \Re^+ \to [r, T]$, $d_{i,j} : \Re^+ \to \{d \in \Re^{k_j} : |d| \leq 1\}$ ($j \neq i$, $i, j = 1, ..., n$) then one can in principle determine from (4.2) the solution $x(t) = (x_1(t), ..., x_n(t))' \in \Re^m$ for $t \in (0, r]$ with $x(\tau) = (x_1(\tau), ..., x_n(\tau))' = x_0(\tau)$ for all $\tau \in [-T, 0]$ using the following procedure:

Step 1:

First determine the solution $x_i(t)$ for $t \in (0, r]$ and for all players with $i \in J_m$ who are using Consistent Backward-looking expectations with respect to the Nash equilibrium point $q^* \in S$. In this case, we are in a position to determine the components of the solution $x_j(t)$ for $j \in J_m$ and $t \in (0, r]$ by means of (4.2). Furthermore, in this case inequality (3.9) holds for $t \in (0, r]$, $i \in J_m$ and consequently there exists function $G_m \in N$ such that:

$$\sup_{t \in [0,r]} |x_i(t)| \leq G_m(\|x\|_{[-T,0]}), \text{ for } i \in J_m \tag{4.3}$$

Step 2:

Next determine the solution $x_i(t)$ for $t \in (0, r]$ and for all players with $i \in J_{m-1}$ who are using Rational-Consistent Backward-looking expectations, $q_{i,j}^{\exp}(t) \in S_j$ if $j \in J_m$ and Consistent Backward-looking expectations, $q_{i,j}^{\exp}(t) \in S_j$ if otherwise. In this case, (3.8) implies that the following inequality holds for all $t \in (0, r]$:

$$|x_i(t)| \leq \max_{\substack{j \neq i \\ j \notin J_m}} \widetilde{\gamma}_{i,j}\left(\|x_j\|_{[t-T, t-r]}\right) + \max_{\substack{j \neq i \\ j \in J_m}} \widetilde{\gamma}_{i,j}\left(\|x_j\|_{[t-T, t]}\right) \tag{4.4}$$

However, the components of the solution $x_j(t)$ for $j \in J_m$ and $t \in (0, r]$ have been determined by Step 1. Therefore, we are in a position to determine the components of the solution $x_j(t)$ for $j \in J_{m-1}$ and $t \in (0, r]$ by means of (4.2). Using (4.3) and (4.4) we obtain the existence of a function $G_{m-1} \in N$ such that:

$$\sup_{t \in [0,r]} |x_i(t)| \leq G_{m-1}(\|x\|_{[-T,0]}), \text{ for } i \in J_m \cup J_{m-1} \tag{4.5}$$

Step $k$ ($3 \leq k \leq m$):

We determine the solution $x_i(t)$ for $t \in (0, r]$ and for all players with $i \in J_{m+1-k}$ who are using Rational-Consistent Backward-looking expectations, $q_{i,j}^{\exp}(t) \in S_j$ if $j \in J_l$ with $l > m+1-k$ and Consistent Backward-looking expectations, $q_{i,j}^{\exp}(t) \in S_j$ if otherwise. In this case, (3.8) implies that the following inequality holds for all $t \in (0, r]$:

$$|x_i(t)| \leq \max_{\substack{j \neq i \\ j \in J_l \\ l \leq m+1-k}} \widetilde{\gamma}_{i,j}\left(\|x_j\|_{[t-T, t-r]}\right) + \max_{\substack{j \neq i \\ j \in J_l \\ l > m+1-k}} \widetilde{\gamma}_{i,j}\left(\|x_j\|_{[t-T, t]}\right) \tag{4.6}$$

However, the components of the solution $x_j(t)$ for $j \in J_l$ with $l > m+1-k$ and $t \in (0, r]$ have been determined by previous steps. Therefore, we are in a position to determine the components of the solution $x_j(t)$ for $j \in J_{m+1-k}$ and $t \in (0, r]$ by means of (4.2). Moreover, by virtue of previous steps there exists there exists function $G_{m+2-k} \in N$ such that:

$$\sup_{t \in [0,r]} |x_i(t)| \leq G_{m+2-k}(\|x\|_{[-T,0]}), \text{ for } i \in \bigcup_{l=m+2-k, ..., m} J_l \tag{4.7}$$

Using (4.6) and (4.7) we obtain the existence of a function $G_{m+1-k} \in \mathcal{N}$ such that:

$$\sup_{t \in [0,r]} |x_i(t)| \leq G_{m+1-k}\left(\|x\|_{[-T,0]}\right), \text{ for } i \in \bigcup_{l=m+1-k,\ldots,m} J_l \quad (4.8)$$

After the completion of the $m$ steps we have determined all components of the solution $x_j(t)$ for $j = 1,\ldots,n$ and $t \in (0,r]$. Moreover, we have also constructed a function $G \in \mathcal{N}$ such that:

$$\sup_{t \in [0,r]} |x(t)| \leq G\left(\|x\|_{[-T,0]}\right) \quad (4.9)$$

Without loss of generality we may assume that $G(s) \geq s$ for all $s \geq 0$. Moreover, by using (4.9) we may conclude exactly as in the previous section that estimates (3.11), (3.12) and (3.13) hold.

The proof of the following theorem is exactly the same with the proof of the Theorem 3.1 and therefore is omitted.

**Theorem 4.2:** $0 \in \mathcal{X}$ is Robustly Globally Asymptotically Stable for system (4.2) under hypothesis (H'), if there exists $\omega > 1$ such that the set of conditions (3.15) holds for each $p = 2,\ldots,n$ and for all $i_j \in \{1,\ldots,n\}$, $i_j \neq i_k$ if $j \neq k$, where $\gamma_{i,j}(s) := \omega \tilde{\gamma}_{i,j}(\omega s)$.

It should be emphasized that the parameters $T \geq r > 0$, which are involved in the definition of the Consistent Backward-looking expectation (Definition 2.1), play no role in the small-gain conditions. Moreover, the number $m$ of the index sets $J_l \subseteq \{1,\ldots,n\}$ involved in hypothesis (H') or the particular members of each index set play absolutely no role in the small-gain conditions (3.15). Furthermore, all these parameters are allowed to change with time: there is no need to assume that these parameters remain constant. Consequently, the small-gain conditions can help us to decide whether the Nash equilibrium point is robustly stable **without any knowledge of the expectation rules**. Again, the small-gain conditions (3.15) demand knowledge of the Nash equilibrium point $q^* \in S$ and the best reply mappings $f_i : S_1 \times \ldots \times S_{i-1} \times S_{i+1} \times \ldots \times S_n \to S_i$ for $1 < i < n$, $n > 3$ and $f_1 : S_2 \times \ldots \times S_n \to S_1$, $f_n : S_1 \times \ldots \times S_{n-1} \to S_n$ for which inequalities (3.8) hold.

Finally, it should be noted that for a specific strategic game, even less demanding hypotheses than hypothesis (H') can be used in order to guarantee that system (3.3) gives an infinite-dimensional dynamical system with state space $\mathcal{X}$ being the normed linear space of bounded functions $x : [-T,0] \to \Re^N$, where $N = k_1 + \ldots + k_n$ with norm $\|x\|_{\mathcal{X}} = \sup_{-T \leq \tau \leq 0} |x(\tau)|$. This can be done by exploiting special properties of the best reply mappings $f_i : S_1 \times \ldots \times S_{i-1} \times S_{i+1} \times \ldots \times S_n \to S_i$ for $1 < i < n$, $n > 3$ and $f_1 : S_2 \times \ldots \times S_n \to S_1$, $f_n : S_1 \times \ldots \times S_{n-1} \to S_n$ (e.g., if some of the functions are independent of certain arguments).

## 5. Conclusions

In this work, advanced stability methods have been used in order to provide sufficient conditions, called cyclic small gain conditions, which guarantee robust global asymptotic stability of the Nash equilibrium in dynamic games. The obtained results are powerful because they can be applied to uncertain models for which the players form consistent expectations based on the history of the game. In addition, by formulating dynamic game-theoretical models by means of Functional Difference Equations, it is possible to obtain all features of continuous-time and discrete-time models. A Cournot oligopoly game has been used in order to illustrate the theoretical results.

Future research can address the economic meaning of small-gain conditions to other games used in economic research (e.g., the study of the stability properties of the Walrasian equilibrium of an abstract economy). A step towards this research direction is the fact that the results presented in this work can be directly extended to the case where the best reply mappings are set-valued maps instead of functions, i.e., $f_i(q_{-i}) \subseteq S_i$. However, in this case inequalities (3.8) must be modified in the following way:

$$\left| p - q_i^* \right| \leq \max_{j \neq i} \widetilde{\gamma}_{i,j}\left(\left| q_j - q_j^* \right|\right), \text{ for all } p \in f_i(q_{-i}) \text{ and } i = 1,...,n \qquad (3.8')$$

The above set of inequalities directly implies that the Nash equilibrium point satisfies $f_i(q_{-i}^*) = \{q_i^*\}$ for all $i = 1,...,n$. With this modification, Theorem 3.1, Corollary 3.2, Corollary 3.3 and Theorem 4.2 hold in this case as well. Future research can also address the issue of studying dynamic game-theoretical models by means of recent results in hybrid systems theory (see [24,25]) or the issue of stabilization of Nash equilibria for dynamic game-theoretical models by means of nonlinear feedback laws, using recently proposed methodologies (see for example [22,23]).

**Acknowledgments**

The work of Z.-P. Jiang has been supported in part by NSF grants ECS-0093176 and DMS-0906659.

# References

[1] Agiza, H. N., G. I. Bischi and M. Kopel, "Multistability in a Dynamic Cournot Game with Three Oligopolists", *Mathematics and Computers in Simulation*, 51, 1999, 63-90.
[2] Athanasiou, G., I. Karafyllis and S. Kotsios, "Price Stabilization Using Buffer Stocks", *Journal of Economic Dynamics and Control*, 32(4), 2008, 1212-1235.
[3] Bischi, G. and M. Kopel, "Equilibrium selection in a nonlinear duopoly game with adaptive expectations", *Journal of Economic Behavior & Organization*, 46, 2001, 73-100.
[4] Chiarella, C. and F. Szidarovszky, "Dynamic oligopolies without full information and with continuously distributed time lags", *Journal of Economic Behavior & Organization*, 54, 2004, 495-511.
[5] Droste, E., C. Hommes, and J. Tuinstra, "Endogenous fluctuations under evolutionary pressure in Cournot competition", *Games and Economic Behavior*, 40, 2002, 232–269.
[6] Forgo, F., J. Szep and F. Szidarovszky, *Introduction to the Theory of Games. Concepts, Methods, Applications*, Kluwer Academic Publishers, 1999.
[7] Grune, L., "Input-to-State Dynamical Stability and its Lyapunov Function Characterization", *IEEE Transactions on Automatic Control*, 47(9), 2002, 1499-1504.
[8] Guo, Q. and P. Niu, "Some Theorems on Existence and Uniqueness of Fixed Points for Decreasing Operators", *Computers and Mathematics with Applications*, 57, 2009, 1515-1521.
[9] Hale, J. K. and S. M. V. Lunel, *Introduction to Functional Differential Equations*, Springer-Verlag, New York, 1993.
[10] Hommes, C. H., "On the consistency of backward-looking expectations: The case of the cobweb", *Journal of Economic Behavior & Organization*, 33, 1998, 333-362.
[11] Huang, W., "Information lag and dynamic stability", *Journal of Mathematical Economics*, 44, 2008, 513-529.
[12] Ito, H. and Z.-P. Jiang, "Nonlinear Small-Gain Condition Covering iISS Systems: Necessity and Sufficiency from a Lyapunov Perspective", *Proceedings of the 44th IEEE Conference on Decision and Control and European Control Conference 2005*, Seville, Spain, 2005, 355-360.
[13] Ito, H. and Z.-P. Jiang, "Necessary and Sufficient Small-Gain Conditions for Integral Input-to-State Stable Systems: A Lyapunov Perspective", *IEEE Transactions on Automatic Control*, 54(10), 2009, 2389-2404.
[14] Jiang, Z.P., A. Teel and L. Praly, "Small-Gain Theorem for ISS Systems and Applications", *Mathematics of Control, Signals and Systems*, 7, 1994, 95-120.
[15] Karafyllis, I., "A System-Theoretic Framework for a Wide Class of Systems I: Applications to Numerical Analysis", *Journal of Mathematical Analysis and Applications*, 328(2), 2007, 876-899.
[16] Karafyllis, I., "A System-Theoretic Framework for a Wide Class of Systems II: Input-to-Output Stability", *Journal of Mathematical Analysis and Applications*, 328(1), 2007, 466-486.
[17] Karafyllis, I. and Z.-P. Jiang, "A Small-Gain Theorem for a Wide Class of Feedback Systems with Control Applications", *SIAM Journal on Control and Optimization*, 46(4), 2007, 1483-1517.
[18] Karafyllis, I., P. Pepe and Z.-P. Jiang, "Stability Results for Systems Described by Coupled Retarded Functional Differential Equations and Functional Difference Equations", *Nonlinear Analysis, Theory, Methods and Applications*, 71(7-8), 2009, 3339-3362.
[19] Karafyllis, I. and Z.-P. Jiang, "A Vector Small-Gain Theorem for General Nonlinear Control Systems", *Proceedings of the 48th IEEE Conference on Decision and Control 2009*, Shanghai, China, 2009, 7996-8001. Also available in http://arxiv.org/abs/0904.0755.
[20] Lin, G., and Y. Hong, "Delay induced oscillation in predator-prey system with Beddington-DeAngelis functional response", *Applied Mathematics & Computation*, 190, 2007, 1297-1311.

[21] Luenberger, D. G., "Complete Stability of Noncooperative Games", *Journal of Optimization Theory and Applications*, 25(4), 1978, 485-505.
[22] Mazenc, F. and P.-A. Bliman, "Backstepping Design for Time-Delay Nonlinear Systems", *IEEE Transactions on Automatic Control*, 51, 2006, 149-154.
[23] Mazenc, F., M. Malisoff and Z. Lin, "Further Results on Input-to-State Stability for Nonlinear Systems with Delayed Feedbacks", *Automatica*, 44, 2008, 2415-2421.
[24] Nesic, D. and A. R. Teel, "Input output stability properties of networked control systems", *IEEE Transactions on Automatic Control*, 49, 2004, 1650-1667.
[25] Nesic, D. and A. R. Teel, "Input to state stability of networked control systems", *Automatica*, 40, 2004, 2121-2128.
[26] Sontag, E.D., "Smooth Stabilization Implies Coprime Factorization", *IEEE Transactions on Automatic Control*, 34, 1989, 435-443.
[27] Sontag, E.D. and Y. Wang, "Notions of Input to Output Stability", *Systems and Control Letters*, 38, 1999, 235-248.
[28] Sontag, E.D. and Y. Wang, "Lyapunov Characterizations of Input-to-Output Stability", *SIAM Journal on Control and Optimization*, 39, 2001, 226-249.
[29] Panchuk, A. and T. Puu, "Stability in a Non-Autonomous Iterative System: An Application to Oligopoly", *Computers and Mathematics with Applications*, 58, 2009, 2022-2034.
[30] Ok, E. A., *Real Analysis with Economic Applications*, Princeton University Press, 2007.
[31] Okuguchi, K., and T. Yamazaki, "Global stability of unique Nash equilibrium in Cournot oligopoly and rent-seeking game", *Journal of Economic Dynamics and Control*, 32, 2008, 1204-1211.
[32] Pepe, P., "The Liapunov's second method for continuous time difference equations", *International Journal of Robust and Nonlinear Control*, 13, 2003, 1389-1405.
[33] Szidarovszky, F. and W. Li, "A Note on the Stability of a Cournot-Nash Equilibrium: the Multiproduct Case with Adaptive Expectations", *Journal of Mathematical Economics*, 33, 2000, 101-107.

# Appendix

**Proof of Theorem 2.2:** Notice that (2.10) implies that the followings equations hold for all $t \geq 0$:

$$q_i(t) = \theta_i(t) \Pr_{[0,Q_i]}\left(q_i(t-\tau_i(t))\right) + (1-\theta_i(t)) \Pr_{[0,Q_i]}\left(\frac{ab-c_i}{2b+K_i} - \frac{b}{2b+K_i}\sum_{j \neq i} q_{i,j}^{\exp}(t)\right), \quad i = 1,\ldots,n \quad (A.1)$$

$$q_{i,j}^{\exp}(t) = \Pr_{[0,Q_j]}\left(q_j^* + d_{i,j}(t)\|q_j - q_j^*\|_{[t-T,t-r]}\right)$$

Using the fact that $\left|\Pr_U(x) - \Pr_U(y)\right| \leq |x-y|$ for every $x, y \in \Re^n$, where $U \subseteq \Re^n$ is a closed convex set, in conjunction with $q_i^* = \Pr_{[0,Q_i]}\left(\frac{ab-c_i}{2b+K_i} - \frac{b}{2b+K_i}\sum_{j \neq i} q_j^*\right)$ and $\theta_i(t) \in [0,\Theta]$ with $\Theta < 1$, we obtain from (A.1) for $i = 1,\ldots,n$:

$$\left|q_i(t) - q_i^*\right| \leq \theta_i(t)\left|q_i(t-\tau_i(t)) - q_i^*\right| + (1-\theta_i(t))\left|\Pr_{[0,Q_i]}\left(\frac{ab-c_i}{2b+K_i} - \frac{b}{2b+K_i}\sum_{j \neq i} q_{i,j}^{\exp}(t)\right) - \Pr_{[0,Q_i]}\left(\frac{ab-c_i}{2b+K_i} - \frac{b}{2b+K_i}\sum_{j \neq i} q_j^*\right)\right| \leq$$

$$\leq \theta_i(t)\left|q_i(t-\tau_i(t)) - q_i^*\right| + \frac{(1-\theta_i(t))b}{2b+K_i}\left|\sum_{j \neq i}(q_{i,j}^{\exp}(t) - q_j^*)\right| \leq \theta_i(t)\left|q_i(t-\tau_i(t)) - q_i^*\right| + \frac{(1-\theta_i(t))b}{2b+K_i}\sum_{j \neq i}\left|q_{i,j}^{\exp}(t) - q_j^*\right|$$

(A.2)

Combining (A.2) with (2.8) and using definitions $R_i = \frac{b}{2b+K_i} > 0$, $x_i(t) = \frac{q_i(t) - q_i^*}{Q_i}$, $g_{i,j} = \frac{Q_j}{Q_i} > 0$ ($j \neq i$, $i = 1,\ldots,n$), we conclude that for every $x_0 \in X$ and for every set of functions $\theta_i : \Re^+ \to [0,\Theta]$, $\tau_i : \Re^+ \to [r,T]$, $d_{i,j} : \Re^+ \to [-1,1]$ ( $j \neq i$, $i,j = 1,...,n$ ), the solution $x_t \in X$ of (2.11) with initial condition $x_{t_0} = x_0 \in X$ corresponding to inputs $\theta_i : \Re^+ \to [0, \Theta]$, $\tau_i : \Re^+ \to [r,T]$, $d_{i,j} : \Re^+ \to [-1,1]$ ( $j \neq i$, $i,j = 1,...,n$ ) satisfies:

$$|x_i(t)| \leq \theta_i(t) \|x_i\|_{[t-T,t-r]} + (1-\theta_i(t)) R_i \sum_{j \neq i} g_{i,j} \|x_j\|_{[t-T,t-r]} \quad \text{for all } t \geq t_0 \tag{A.3}$$

We notice that system (2.11) is an autonomous uncertain dynamical system in the sense described in [15,16,17]. Next we show that $0 \in X$ is a robust equilibrium point for system (2.11) in the sense described in [15,16,17], i.e., for every $\varepsilon > 0$, $T > 0$ there exists $\delta := \delta(\varepsilon, T) > 0$ such that if $\|x_0\|_X \leq \delta$ then for every set of functions $\theta_i : \Re^+ \to [0, \Theta]$, $\tau_i : \Re^+ \to [r,T]$, $d_{i,j} : \Re^+ \to [-1,1]$ ( $j \neq i$, $i,j = 1,...,n$ ), the solution $x_t \in X$ of (2.11) with initial condition $x_0 \in X$ corresponding to inputs $\theta_i : \Re^+ \to [0, \Theta]$, $\tau_i : \Re^+ \to [r,T]$, $d_{i,j} : \Re^+ \to [-1,1]$ ( $j \neq i$, $i,j = 1,...,n$ ) satisfies $\|x_t\|_X = \|x\|_{[t-T,t]} \leq \varepsilon$ for all $t \in [0,T]$. To see this, notice that (A.3) implies the existence of a constant $G \geq 0$ such that:

$$\sup_{t \in [0,r]} |x(t)| = \|x\|_{[0,r]} \leq G \|x\|_{[-T,0]} \tag{A.4}$$

Without loss of generality we may assume that $G \geq 1$. Inequality (A.4) implies that

$$\|x\|_{[t-T,t]} \leq G \|x\|_{[-T,0]}, \text{ for all } t \in [0,r] \tag{A.5}$$

Working in this way and using induction we may establish that for every positive integer $k > 0$ it holds that

$$\|x\|_{[t-T,t]} \leq G^k \|x\|_{[-T,0]}, \text{ for all } t \in [0, kr] \tag{A.6}$$

Therefore (A.6) implies that

$$\|x_t\|_X = \|x\|_{[t-T,t]} \leq G^{1+[t/r]} \|x\|_{[-T,0]} = G^{1+[t/r]} \|x_0\|_X, \text{ for all } t \geq 0 \tag{A.7}$$

where $[t/r]$ denotes the integer part of $t/r$. Consequently, (A.7) implies that for every $\varepsilon > 0$, $T > 0$ there exists $\delta := \delta(\varepsilon, T) = \varepsilon G^{-1-[T/r]} > 0$ such that if $\|x_0\|_X \leq \delta$ then for every set of functions $\theta_i : \Re^+ \to [0, \Theta]$, $\tau_i : \Re^+ \to [r,T]$, $d_{i,j} : \Re^+ \to [-1,1]$ ( $j \neq i$, $i,j = 1,...,n$ ), the solution $x_t \in X$ of (2.11) with initial condition $x_0 \in X$ corresponding to inputs $\theta_i : \Re^+ \to [0, \Theta]$, $\tau_i : \Re^+ \to [r,T]$, $d_{i,j} : \Re^+ \to [-1,1]$ ( $j \neq i$, $i,j = 1,...,n$ ) satisfies $\|x_t\|_X = \|x\|_{[t-T,t]} \leq \varepsilon$ for all $t \in [0,T]$. Therefore $0 \in X$ is a robust equilibrium point for system (2.11) in the sense described in [15,16,17].

The reader should notice that inequality (A.3) implies the following inequality for all $i = 1,...,n$ and $\mu > \Theta$:

$$|x_i(t)| \leq \max\left\{ \mu \|x_i\|_{[t-T,t-r]}, \frac{\mu - \mu \Theta}{\mu - \Theta} R_i \sum_{j \neq i} g_{i,j} \|x_j\|_{[t-T,t-r]} \right\} \text{ for all } t \geq t_0 \tag{A.8}$$

Let $\sigma > 0$ and consider the family of functionals $V_i : X \to \Re^+$, $i = 1,...,n$ defined by:

$$V_i(x) = \sup_{-T \leq \tau \leq 0} Q_i |x_i(\tau)| \exp(\sigma \tau) \tag{A.9}$$

Let $h \in (0, r)$ and $t \geq 0$ be arbitrary. Definition (A.9) and inequality (A.8) imply that:

$$V_i(x_{t+h}) = \sup_{-T \leq \tau \leq 0} Q_i |x_i(t+h+\tau)| \exp(\sigma\tau) =$$

$$= \sup_{t+h-T \leq s \leq t+h} Q_i |x_i(s)| \exp(\sigma(s-t-h))$$

$$\leq \max \left\{ \sup_{t+h-T \leq s \leq t} Q_i |x_i(s)| \exp(\sigma(s-t-h)), \sup_{t \leq s \leq t+h} Q_i |x_i(s)| \exp(\sigma(s-t-h)) \right\}$$

$$\leq \max \left\{ \begin{array}{l} \exp(-\sigma h) V_i(x_t), \mu Q_i \sup_{t \leq s \leq t+h} \|x_i\|_{[s-T, s-r]} \exp(\sigma(s-t-h)), \\ \dfrac{\mu - \mu\Theta}{\mu - \Theta} R_i \sup_{t \leq s \leq t+h} \sum_{j \neq i} Q_i g_{i,j} \|x_j\|_{[s-T, s-r]} \exp(\sigma(s-t-h)) \end{array} \right\}$$

Using definition (A.9) and the facts that $g_{i,j} = \dfrac{Q_j}{Q_i} > 0$ ($j \neq i$, $i, j = 1, ..., n$) and $\|x_j\|_{[s-T, s-r]} = \sup_{s-T \leq \tau \leq s-r} |x_j(\tau)| = \sup_{-T \leq w \leq -r} |x_j(s+w)|$, we obtain from the above inequality:

$$V_i(x_{t+h}) \leq$$

$$\leq \max \left\{ \begin{array}{l} \exp(-\sigma h) V_i(x_t), \mu \sup_{t \leq s \leq t+h} \sup_{-T \leq w \leq -r} Q_i |x_i(s+w)| \exp(\sigma w) \exp(-\sigma w) \exp(\sigma(s-t-h)), \\ \dfrac{\mu - \mu\Theta}{\mu - \Theta} R_i \sup_{t \leq s \leq t+h} \sum_{j \neq i} Q_j \sup_{-T \leq w \leq -r} |x_j(s+w)| \exp(\sigma w) \exp(-\sigma w) \exp(\sigma(s-t-h)) \end{array} \right\}$$

$$\leq \max \left\{ \exp(-\sigma h) V_i(x_t), \mu \exp(\sigma T) \sup_{t \leq s \leq t+h} V_i(x_s), \dfrac{\mu - \mu\Theta}{\mu - \Theta} R_i \exp(\sigma T) \sup_{t \leq s \leq t+h} \sum_{j \neq i} V_j(x_s) \right\}$$

$$\leq \max \left\{ \exp(-\sigma h) V_i(x_t), \mu \exp(\sigma T) \sup_{t \leq s \leq t+h} V_i(x_s), \dfrac{\mu - \mu\Theta}{\mu - \Theta} R_i (n-1) \exp(\sigma T) \max_{j \neq i} \sup_{t \leq s \leq t+h} V_j(x_s) \right\}$$

Consequently, for every $i = 1, ..., n$, $\mu > \Theta$, $\sigma > 0$, $h \in (0, r)$ and $t \geq 0$ it holds that:

$$V_i(x_{t+h}) \leq \max \left\{ \exp(-\sigma h) V_i(x_t), \mu \exp(\sigma T) \sup_{t \leq s \leq t+h} V_i(x_s), \dfrac{\mu - \mu\Theta}{\mu - \Theta} R_i (n-1) \exp(\sigma T) \max_{j \neq i} \sup_{t \leq s \leq t+h} V_j(x_s) \right\} \quad \text{(A.10)}$$

Using induction and (A.10), we can show that for every $i = 1, ..., n$, $\mu > \Theta$, $\sigma > 0$, $h \in (0, r)$, $t \geq 0$ and for every non-negative integer $k \geq 0$, it holds that:

$$V_i(x_{t+kh}) \leq \max \left\{ \exp(-\sigma k h) V_i(x_t), \mu \exp(\sigma T) \sup_{t \leq s \leq t+kh} V_i(x_s), \dfrac{\mu - \mu\Theta}{\mu - \Theta} R_i (n-1) \exp(\sigma T) \max_{j \neq i} \sup_{t \leq s \leq t+kh} V_j(x_s) \right\}$$
(A.11)

Therefore, (A.11) implies that for every $i = 1, ..., n$, $\mu > \Theta$, $\sigma > 0$ and $t \geq 0$ the following inequality holds:

$$V_i(x_t) \leq \max \left\{ \exp(-\sigma t) V_i(x_0), \mu \exp(\sigma T) \sup_{0 \leq s \leq t} V_i(x_s), \dfrac{\mu - \mu\Theta}{\mu - \Theta} R_i (n-1) \exp(\sigma T) \max_{j \neq i} \sup_{0 \leq s \leq t} V_j(x_s) \right\} \quad \text{(A.12)}$$

Next, we assume that $\sigma < T^{-1} \ln(2)$. The reader should notice that definition (A.8) implies that $\dfrac{2}{Q_i} V_i(x_t) \geq \sup_{-T \leq \tau \leq 0} |x_i(t+\tau)|$ for $\sigma < T^{-1} \ln(2)$ and consequently:

$$\|x_t\|_X = \|x\|_{[t-T,t]} \le \sum_{i=1}^n \sup_{-T \le \tau \le 0} |x_i(t+\tau)| \le \frac{2}{q} \sum_{i=1}^n V_i(x_t) \tag{A.13}$$

where $q := \min_{i=1,\ldots,n} Q_i$. It follows from (A.12), (A.13) and definition (A.9) (which implies $V_i(x_t) \le Q\|x_t\|_X$ for $Q := \max_{i=1,\ldots,n} Q_i$) that the following inequalities hold for every $i = 1,\ldots,n$, $\mu > \Theta$, $\sigma > 0$ and $t \ge 0$ with $\sigma < T^{-1} \ln(2)$:

$$V_i(x_t) \le \max\left\{ Q\exp(-\sigma t)L(x_0), \mu\exp(\sigma T) \sup_{0 \le s \le t} V_i(x_s), \frac{\mu - \mu\Theta}{\mu - \Theta} R_i(n-1)\exp(\sigma T) \max_{j \ne i} \sup_{0 \le s \le t} V_j(x_s) \right\} \tag{A.14}$$

$$L(x_t) \le \max\left\{ \frac{4nQ}{q} \|x_0\|_X, \frac{4n}{q}(R(n-1)+1) \sum_{j=1}^n \sup_{0 \le s \le t} V_j(x_s) \right\} \tag{A.15}$$

where $L(x) := \|x\|_X$ and $R := \max_{i=1,\ldots,n} R_i$. It follows from (A.13), (A.14), (A.15) and Theorem 3.1 in [19] that $0 \in X$ is Robustly Globally Asymptotically Stable for system (2.11), provided that the following set of conditions holds for each $p = 2,\ldots,n$:

$$R_{i_1} \ldots R_{i_p}(n-1)^p \left(\frac{\mu - \mu\Theta}{\mu - \Theta}\right)^p \exp(p\sigma T) < 1 \tag{A.16}$$

for all $i_j \in \{1,\ldots,n\}$, $i_j \ne i_k$ if $j \ne k$ and

$$\mu \exp(\sigma T) < 1 \tag{A.17}$$

Notice that if conditions (2.15) hold for each $p = 2,\ldots,n$ and for all $i_j \in \{1,\ldots,n\}$, $i_j \ne i_k$ if $j \ne k$ then the conditions (A.16), (A.17) hold for each $p = 2,\ldots,n$ and for all $i_j \in \{1,\ldots,n\}$, $i_j \ne i_k$ if $j \ne k$ for sufficiently small $\sigma > 0$ and for $\mu \in (\Theta,1)$ sufficiently close to 1. The proof is complete. ◁

**Proof of Theorem 3.1:** We first notice that system (3.7) is an autonomous dynamical system in the sense described in [15,16,17]. Next we show that $0 \in X$ is a robust equilibrium point for system (3.7) in the sense described in [15,16,17], i.e., for every $\varepsilon > 0$, $T > 0$ there exists $\delta := \delta(\varepsilon,T) > 0$ such that if $\|x_0\|_X \le \delta$ then for every set of functions $\theta_i : \Re^+ \to [0,\Theta]$, $\tau_i : \Re^+ \to [r,T]$, $d_{i,j} : \Re^+ \to \{d \in \Re^{k_j} : |d| \le 1\}$ ($j \ne i$, $i,j = 1,\ldots,n$), the solution $x_t \in X$ of (3.7) with initial condition $x_0 \in X$ corresponding to inputs $\theta_i : \Re^+ \to [0,\Theta]$, $\tau_i : \Re^+ \to [r,T]$, $d_{i,j} : \Re^+ \to \{d \in \Re^{k_j} : |d| \le 1\}$ ($j \ne i$, $i,j = 1,\ldots,n$) satisfies $\|x_t\|_X = \|x\|_{[t-T,t]} \le \varepsilon$ for all $t \in [0,T]$. Without loss of generality we may assume that the function $G \in N$ involved in (3.13) is a strictly increasing function. Define $\kappa(s) := G^{(1+[T/r])}(s)$ for every $T > 0$, which is a strictly increasing, continuous function with $\kappa(0) = 0$ and $\lim_{s \to +\infty} \kappa(s) = +\infty$ (recall that $G(s) \ge s$ for all $s \ge 0$) and define $\kappa^{-1} : \Re^+ \to \Re^+$ to be the inverse function of $\kappa$ on $\Re^+$. Indeed, (3.13) implies that for every $\varepsilon > 0$, $T > 0$ there exists $\delta := \delta(\varepsilon,T) = \kappa^{-1}(\varepsilon) > 0$ such that if $\|x_0\|_X \le \delta$ then for every set of functions $\theta_i : \Re^+ \to [0,\Theta]$, $\tau_i : \Re^+ \to [r,T]$, $d_{i,j} : \Re^+ \to \{d \in \Re^{k_j} : |d| \le 1\}$ ($j \ne i$, $i,j = 1,\ldots,n$), the solution $x_t \in X$ of (3.7) with initial condition $x_0 \in X$ corresponding to inputs $\theta_i : \Re^+ \to [0,\Theta]$, $\tau_i : \Re^+ \to [r,T]$, $d_{i,j} : \Re^+ \to \{d \in \Re^{k_j} : |d| \le 1\}$ ($j \ne i$, $i,j = 1,\ldots,n$) satisfies $\|x_t\|_X = \|x\|_{[t-T,t]} \le \varepsilon$ for all $t \in [0,T]$. Therefore $0 \in X$ is a robust equilibrium point for system (3.7) in the sense described in [15,16,17].

Let $\sigma > 0$ and consider the family of functionals $V_i : X \to \Re^+$, $i = 1,\ldots,n$ defined by:

$$V_i(x) = \sup_{-T \leq \tau \leq 0} |x_i(\tau)| \exp(\sigma\tau) \qquad (A.18)$$

Using definition (A.18) and (3.9) we obtain for every $h \in (0,T)$, $i = 1,...,n$, $\mu > \Theta$ and $t \geq 0$:

$$\begin{aligned}
V_i(x_{t+h}) &= \sup_{-T \leq \tau \leq 0} |x_i(t+h+\tau)| \exp(\sigma\tau) = \\
&= \sup_{t+h-T \leq s \leq t+h} |x_i(s)| \exp(\sigma(s-t-h)) \\
&\leq \max\left\{\sup_{t+h-T \leq s \leq t} |x_i(s)| \exp(\sigma(s-t-h)), \sup_{t \leq s \leq t+h} |x_i(s)| \exp(\sigma(s-t-h))\right\} \\
&\leq \max\left\{\exp(-\sigma h) V_i(x_t), \sup_{t \leq s \leq t+h} \mu \|x_i\|_{[s-T,s]} \exp(\sigma(s-t-h)), \sup_{t \leq s \leq t+h} \max_{j \neq i} \frac{\mu - \mu\Theta}{\mu - \Theta} \tilde{\gamma}_{i,j}\left(\|x_j\|_{[s-T,s]}\right) \exp(\sigma(s-t-h))\right\} \\
&= \max\left\{\begin{array}{l} \exp(-\sigma h) V_i(x_t), \sup_{t \leq s \leq t+h} \mu \sup_{-T \leq w \leq 0} |x_i(w+s)| \exp(\sigma w) \exp(-\sigma w) \exp(\sigma(s-t-h)), \\ \sup_{t \leq s \leq t+h} \max_{j \neq i} \frac{\mu - \mu\Theta}{\mu - \Theta} \tilde{\gamma}_{i,j}\left(\sup_{-T \leq w \leq 0} |x_j(w+s)| \exp(\sigma w) \exp(-\sigma w)\right) \exp(\sigma(s-t-h)) \end{array}\right\} \\
&\leq \max\left\{\exp(-\sigma h) V_i(x_t), \sup_{t \leq s \leq t+h} \mu \exp(\sigma T) V_i(x_s), \sup_{t \leq s \leq t+h} \max_{j \neq i} \frac{\mu - \mu\Theta}{\mu - \Theta} \tilde{\gamma}_{i,j}\left(\exp(\sigma T) V_j(x_s)\right)\right\} \\
&\leq \max\left\{\exp(-\sigma h) V_i(x_t), \mu \exp(\sigma T) \sup_{t \leq s \leq t+h} V_i(x_s), \max_{j \neq i} \frac{\mu - \mu\Theta}{\mu - \Theta} \tilde{\gamma}_{i,j}\left(\exp(\sigma T) \sup_{t \leq s \leq t+h} V_j(x_s)\right)\right\}
\end{aligned}$$

Consequently, for every $i = 1,...,n$, $\sigma > 0$, $h \in (0,T)$, $\mu > \Theta$ and $t \geq 0$ it holds that:

$$V_i(x_{t+h}) \leq \max\left\{\exp(-\sigma h) V_i(x_t), \mu \exp(\sigma T) \sup_{t \leq s \leq t+h} V_i(x_s), \max_{j \neq i} \frac{\mu - \mu\Theta}{\mu - \Theta} \tilde{\gamma}_{i,j}\left(\exp(\sigma T) \sup_{t \leq s \leq t+h} V_j(x_s)\right)\right\} \qquad (A.19)$$

Using induction and (A.19), we can show that for every $i = 1,...,n$, $\sigma > 0$, $h \in (0,T)$, $\mu > \Theta$, $t \geq 0$ and for every non-negative integer $k \geq 0$, it holds that:

$$V_i(x_{t+kh}) \leq \max\left\{\exp(-\sigma k h) V_i(x_t), \mu \exp(\sigma T) \sup_{t \leq s \leq t+kh} V_i(x_s), \max_{j \neq i} \frac{\mu - \mu\Theta}{\mu - \Theta} \tilde{\gamma}_{i,j}\left(\exp(\sigma T) \sup_{t \leq s \leq t+kh} V_j(x_s)\right)\right\} \qquad (A.20)$$

Therefore, (A.20) implies that for every $i = 1,...,n$, $\mu > \Theta$, $\sigma > 0$ and $t \geq 0$ the following inequality holds:

$$V_i(x_t) \leq \max\left\{\exp(-\sigma t) V_i(x_0), \mu \exp(\sigma T) \sup_{0 \leq s \leq t} V_i(x_s), \max_{j \neq i} \frac{\mu - \mu\Theta}{\mu - \Theta} \tilde{\gamma}_{i,j}\left(\exp(\sigma T) \sup_{0 \leq s \leq t} V_j(x_s)\right)\right\} \qquad (A.21)$$

The reader should notice that definition (A.18) implies that $2V_i(x_t) \geq \sup_{-T \leq \tau \leq 0} |x_i(t+\tau)|$ for $\sigma \leq T^{-1} \ln(2)$ and consequently:

$$\|x_t\|_\mathcal{X} = \|x\|_{[t-T,t]} \leq \sum_{i=1}^n \sup_{-T \leq \tau \leq 0} |x_i(t+\tau)| \leq 2 \sum_{i=1}^n V_i(x_t) \qquad (A.22)$$

Without loss of generality we may assume $\Theta > 0$. Define $\mu := \frac{\omega\Theta}{\omega - 1 + \Theta}$ and let the constant $\sigma > 0$ satisfy the inequalities: $\sigma \leq T^{-1} \ln(2)$, $\sigma \leq T^{-1} \ln(\omega)$ and $\sigma < T^{-1} \ln\left(\frac{\omega - 1 + \Theta}{\omega\Theta}\right)$, where $\omega > 1$ is the constant involved in the hypotheses of the theorem. Notice that the hypothesis $\Theta < 1$ and previous definitions imply that $\mu \exp(\sigma T) < 1$, $\mu > \Theta$, $\exp(\sigma T) \le \omega$ and $\dfrac{\mu - \mu \Theta}{\mu - \Theta} \le \omega$. It follows from (A.21), (A.22) and definition (A.18) (which implies $V_i(x_t) \le \|x_t\|_{\mathcal{X}}$) that the following inequalities hold for every $i = 1,\ldots,n$ and $t \ge 0$:

$$V_i(x_t) \le \max\left\{ \exp(-\sigma t) L(x_0), B \sup_{0 \le s \le t} V_i(x_s), \max_{j \ne i} \omega \widetilde{\gamma}_{i,j}\left( \omega \sup_{0 \le s \le t} V_j(x_s) \right) \right\} \quad (A.23)$$

$$L(x_t) \le \max\left\{ 2n\|x_0\|_{\mathcal{X}}, 2\sum_{i=1}^{n} \max_{j \ne i} \omega \widetilde{\gamma}_{i,j}\left( \omega \sup_{0 \le s \le t} V_j(x_s) \right) + 2B \sum_{i=1}^{n} \sup_{0 \le s \le t} V_i(x_s) \right\} \quad (A.24)$$

where $L(x) := \|x\|_{\mathcal{X}}$ and $B := \mu \exp(\sigma T) < 1$. It follows from (A.22), (A.23), (A.24) and Theorem 3.1 in [19] that $0 \in \mathcal{X}$ is Robustly Globally Asymptotically Stable for system (3.7), provided that the set of conditions (3.15) holds for each $p = 2,\ldots,n$ and for all $i_j \in \{1,\ldots,n\}$, $i_j \ne i_k$ if $j \ne k$. The proof is complete. ◁